\documentclass[12pt]{article}
\usepackage{amsfonts}
\newcommand{\beqn}{\begin{eqnarray}}
\newcommand{\eeqn}{\end{eqnarray}}

\newtheorem{proposition}{Proposition}
\newtheorem{lemma}{Lemma}

\newcommand{\proof}{\noindent {\it Proof.\/}\ }
\newcommand{\qed}{{\it Q.E.D.\/} \bigskip\par}
\newcommand{\rd}{\partial}
\newcommand{\dfrac}[2]{ \frac{\displaystyle #1}{\displaystyle #2} }
\newcommand{\ssfrac}[2]%
{ \frac{\scriptscriptstyle #1}{\scriptscriptstyle #2} }

\newcommand{\Tr}{\mathop{{\mathrm{Tr}}}}

\newcommand{\const}{{\mathrm{const.}}}

\newcommand{\bbC}{{\mathbb{C}}}
\newcommand{\bbP}{{\mathbb{P}}}
\newcommand{\bbR}{{\mathbb{R}}}
\newcommand{\bbZ}{{\mathbb{Z}}}
\newcommand{\frakg}{{\mathfrak{g}}}
\newcommand{\frakh}{{\mathfrak{h}}}
\newcommand{\calH}{{\mathcal{H}}}
\newcommand{\gtilde}{\tilde{g}}
\begin{document}

\title{Elliptic Calogero-Moser Systems \\
and \\
Isomonodromic Deformations}
\author{Kanehisa Takasaki\\
{\normalsize Department of Fundamental Sciences, Kyoto University}\\
{\normalsize Yoshida, Sakyo-ku, Kyoto 606-8501, Japan}\\
{\normalsize E-mail: takasaki@yukawa.kyoto-u.ac.jp}}
\date{}
\maketitle

\begin{abstract}
We show that various models of the elliptic Calogero-Moser systems 
are accompanied with an isomonodromic system on a torus.  The 
isomonodromic partner is a non-autonomous Hamiltonian system 
defined by the same Hamiltonian.  The role of the time variable 
is played by the modulus of the base torus.  A suitably chosen 
Lax pair (with an elliptic spectral parameter) of the elliptic 
Calogero-Moser system turns out to give a Lax representation of 
the non-autonomous system as well.  This Lax representation 
ensures that the non-autonomous system describes isomonodromic 
deformations of a linear ordinary differential equation on the 
torus on which the spectral parameter of the Lax pair is defined.  
A particularly interesting example is the ``extended twisted 
$BC_\ell$ model'' recently introduced along with some other models 
by Bordner and Sasaki, who remarked that this system is equivalent 
to Inozemtsev's generalized elliptic Calogero-Moser system.  We 
use the ``root type'' Lax pair developed by Bordner et al. to 
formulate the associated isomonodromic system on the torus. 
\end{abstract}
\begin{flushleft}
KUCP-0133\\
math.QA/9905101
\end{flushleft}
%%%%%%%%%%%%%%%%%%%%%%%%%%%%%%%%%%%%%%%%%%%%%%%%%%%%%%%%%%%%%%%%%%%%
\newpage
\renewcommand{\theequation}{\arabic{section}.\arabic{equation}}
%%%%%%%%%%%%%%%%%%%%%%%%%%%%%%%%%%%%%%%%%%%%%%%%%%%%%%%%%%%%%%%%%%%%
\section{Introduction}
\setcounter{equation}{0}

In 1996, Manin \cite{bib:manin96} proposed a new 
expression of the sixth Painlev\'e equation. 
This is a differential equation of the form 
\beqn
    (2 \pi i)^2 \frac{d^2 q}{d\tau^2} 
    = \sum_{a=0}^3 \alpha_a \wp'(q + \omega_a), 
\eeqn
where $\wp'(u)$ is the derivative of the Weierstrass 
$\wp$ function with primitive periods $1$ and $\tau$, 
\beqn
    \wp(u) = \wp(u \mid 1,\tau) 
    = \frac{1}{u^2} 
      + \sum_{(m,n) \not= (0,0)} 
        \left( \frac{1}{(u + m + n\tau)^2} 
        - \frac{1}{(m + n\tau)^2} \right), 
\eeqn
$\omega_a$ ($a = 0,1,2,3$) are the origin and the three 
half-periods of the torus $E_\tau = \bbC/(\bbZ + \tau\bbZ)$, 
\beqn
    \omega_0 = 0, \quad 
    \omega_1 = \frac{1}{2}, \quad 
    \omega_2 = \frac{1}{2} + \frac{\tau}{2}, \quad 
    \omega_3 = \frac{\tau}{2}, 
\eeqn
and $\alpha_a$ ($a = 0,1,2,3$) are the simple linear 
combinations $(\alpha_0,\alpha_1,\alpha_2,\alpha_3) 
= (\alpha, -\beta, \gamma, 1/2 - \delta)$ of 
the four parameters $\alpha$, $\beta$, $\gamma$ 
and $\beta$ of the sixth Painlev\'e equation
\beqn
    \frac{dy^2}{dx^2} 
    &=& \frac{1}{2}\left( \frac{1}{y} + \frac{1}{y-1} 
          + \frac{1}{y-x} \right) \left(\frac{dy}{dx}\right)^2 
        - \left( \frac{1}{x} + \frac{1}{x-1} 
          + \frac{1}{y-x} \right) \frac{dy}{dx} 
    \nonumber \\
    &&  + \frac{y(y-1)(y-x)}{x^2(x-1)^2} 
          \left( \alpha + \beta \frac{x}{y^2} 
          + \gamma \frac{x-1}{(y-1)^2} 
          + \delta \frac{x(x-1)}{(y-x)^2} \right). 
\eeqn
Manin's equation can be written in the Hamiltonian form 
\beqn
    2 \pi i\frac{dq}{d\tau} = p, \quad 
    2 \pi i\frac{dp}{d\tau} = - \frac{\rd \calH}{\rd q} 
\eeqn
with the Hamiltonian 
\beqn
    \calH = \frac{1}{2}p^2 
      - \sum_{a=0}^3 \alpha_a \wp(q + \omega_a). 
\eeqn
Since the Hamiltonian depends on the modulus $\tau$ 
explicitly, this is a non-autonomous Hamiltonian system.  
In this new framework, Manin reconsidered the affine 
Weyl group symmetries of the sixth Painlev\'e equation 
discovered by Okamoto \cite{bib:okamoto87}, solutions 
for special values of $\alpha$, $\beta$, $\gamma$ and 
$\delta$ constructed by Hitchin \cite{bib:hitchin95}, 
etc.  

Manin's equation reveals an unexpected link between 
the Painlev\'e equation and the elliptic Calogero-Moser 
systems, i.e., the Calogero-Moser systems 
\cite{bib:calogero-moser} with elliptic potentials.  
In order to see this relation, we introduce a new 
variable $t$ and formally replace 
$2 \pi i d/d\tau \to d/dt$ 
in the aforementioned equations.  The outcome are 
the autonomous equation 
\beqn
    \frac{d^2 q}{dt^2} 
    = \sum_{a=0}^3 \alpha_a \wp'(q + \omega_a) 
\eeqn
and its Hamiltonian form 
\beqn
    \frac{dq}{dt} = p, \quad 
    \frac{dp}{dt} = - \frac{\rd \calH}{\rd q}. 
\eeqn
If all $\alpha_n$'s take the same value $-g^2/8$, 
one can use an identity of the $\wp$ function to 
rewrite the above equation as: 
\beqn
    \frac{d^2 q}{dt^2} 
    =  - \frac{g^2}{8} \sum_{a=0}^3 \wp'(q + \omega_a) 
    = - g^2 \wp'(2q). 
\eeqn
This is exactly the two-body elliptic Calogero-Moser 
system; the $\ell$-body elliptic Calogero-Moser system 
($A_{\ell-1}$ model) is defined by the Hamiltonian 
\beqn
    \calH = \frac{1}{2} \sum_{j=1}^\ell p_j^2 
      + \frac{g^2}{2} \sum_{j\not= k} \wp(q_j - q_k). 
\eeqn
As Krichever \cite{bib:krichever80} demonstrated, 
this elliptic Calogero-Moser system is an isospectral 
integrable system with a Lax representation 
\beqn
    \frac{\rd L(z)}{\rd t} = [L(z), M(z)], 
\eeqn
where the Lax pair $L(z)$ and $M(z)$ are matrix-valued 
functions of a spectral parameter $z$ on the torus 
$E_\tau$.  Furthermore, the general case falls into 
Inozemtsev's generalization of the elliptic 
Calogero-Moser system \cite{bib:inozemtsev89} 
defined by the Hamiltonian 
\beqn
    \calH = \frac{1}{2} \sum_{j=1}^\ell p_j^2 
       + \frac{g_m^2}{2} \sum_{\epsilon,\epsilon'=\pm 1}
         \sum_{j\not= k} \wp(\epsilon q_j + \epsilon' q_k) 
       + \frac{1}{2} \sum_{j=1}^\ell \sum_{a=0}^3 
         g_a^2 \wp(q_j + \omega_a). 
\eeqn

Levin and Olshanetsky \cite{bib:levin-olsh97}
developed a geometric formulation of isomonodromic 
systems on a general Riemann surface, and 
characterized Manin's equation as an isomonodromic 
system on the torus $E_\tau$.  Their interpretation 
of isomonodromic deformations is based on the notion 
of the Hitchin systems \cite{bib:hitchin87}.  
According to this interpretation, the coordinates 
$q_j$ of Calogero-Moser particles are identified 
with the moduli of an $SU(\ell)$ flat bundle on 
the torus $E_\tau$, and the $L$-matrix $L(z)$ is 
nothing but the Higgs field on this bundle.  
(Such a link between the elliptic Calogero-Moser 
systems and the Hitchin systems was already pointed 
out before their work by Nekrasov \cite{bib:nekrasov95} 
and Enriquez and Rubtsov \cite{bib:enri-rubt95}.) 
Isomonodromic deformations are special deformations 
of these geometric data as the complex structure 
of the base torus (or, equivalently, the modulus 
$\tau$) varies.  This geometric picture suggests 
a wide range of generalizations of isomonodromic 
deformations (see, e.g., the recent work of Levin 
and Olshanetsky \cite{bib:levin-olsh99}).  

Unfortunately, however, it is only the special case 
with $\alpha_0 = \alpha_1 = \alpha_2 = \alpha_3$ 
that was successfully treated in the formulation 
of Levin and Olshanetsky.  This is simply because 
no suitable Lax representation was available for the 
Inozemtsev system.  Inozemtsev \cite{bib:inozemtsev89} 
presented a Lax representation, but it is not suited 
for that purpose.  

Recently, a new type of Lax pair --- the root type 
Lax pair --- was proposed by Bordner et al. \cite{% 
bib:bordner-etal98a,bib:bordner-etal98b,bib:bordner-etal98c} 
for various models of the elliptic Calogero-Moser 
systems including the Inozemtsev system. This is a 
Lax pair constructed on the basis of an underlying 
root system (e.g., the $A_{\ell-1}$ root system for 
the aforementioned elliptic Calogero-Moser system, 
and the $BC_\ell$ root system for the Inozemtsev 
system).  The construction covers not only the 
ordinary elliptic Calogero-Moser systems (the 
``untwisted models'') but also the ``twisted models'' 
introduced by D'Hoker and Phong \cite{bib:dhoker-phong98} 
and their generalizations (the ``extended twisted models'').  
The Inozemtsev system coincides with the extended 
twisted $BC_\ell$ model in the classification of 
Bordner and Sasaki \cite{bib:bordner-etal98c}.   
In particular, the root type Lax pair for the extended 
twisted $BC_1$ model gives a Lax representation to 
the aforementioned isospectral analogue of Manin's 
equation.  

One of the goals of this paper is to show, using the 
root type Lax pair, that each of these elliptic 
Calogero-Moser systems are accompanied with an 
isomonodromic system on a torus.  The fist step of 
the construction is simply to replace the equations 
of motions 
\beqn
    \frac{dq}{dt} = \{q, \calH\}, \quad 
    \frac{dp}{dt} = \{p, \calH\} 
\eeqn
of the elliptic Calogero-Moser system by the 
non-autonomous system 
\beqn
    2 \pi i \frac{dq}{d\tau} = \{q, \calH\}, \quad 
    2 \pi i \frac{dp}{d\tau} = \{p, \calH\} 
\eeqn
with the same Hamiltonian $\calH$.  We then rewrite 
this non-autonomous system into a Lax equation of 
the form 
\beqn
    2 \pi i\frac{\rd L(z)}{\rd \tau} 
    + \frac{\rd M(z)}{\rd z} 
    = [L(z), M(z)] 
\eeqn
using a root type Lax pair $L(z)$ and $M(z)$. This 
Lax equation implies the Frobenius integrability 
of the linear system 
\beqn
    \frac{\rd Y(z)}{\rd z} = L(z) Y(z), \quad 
    2 \pi i\frac{\rd L(z)}{\rd \tau} + M(z) Y(z) = 0, 
\eeqn
from which one can deduce that the non-autonomous 
system is an isomonodromic system on the torus $E_\tau$.  

Actually, we shall use the root type Lax pair made of 
slightly different building blocks.  The root type 
Lax pairs, like the previously known Lax pairs, contain 
complex analytic functions $x(u,z)$, $y(u,z)$, etc. 
that satisfy special functional equations (called the 
``Calogero functional equations'' \cite{bib:calogero76}).  
Bordner et al. use the Weierstrass sigma function to 
construct those functions.  We use the Jacobi theta 
function $\theta_1$ instead.  This is inspired by 
the work of Levin and Olshanetsky, who used 
substantially the same function to construct the 
$L$-matrix (i.e., the Higgs field in their framework) 
for isomonodromic systems on a torus. This minuscule 
difference is rather crucial for deriving an 
isomonodromic Lax equation as above.  

The functions $x(u,z)$ and $y(u,z)$ that we use are, 
in fact, identical to the functions that Felder and 
Wieczerkowski \cite{bib:feld-wiec94} used in their 
study on the Knizhnik-Zamolodchikov-Bernard (KZB) 
equation \cite{bib:bernard88}.  This is by no means 
a coincidence.  As Levin and Olshanetsky stressed, 
the KZB equation and the Hitchin system (or, rather, 
its isomonodromic version) are closely related.  

In order to illustrate that our method also works 
for some other cases, we show a construction of an 
isomonodromic analogue for the ``spin generalization'' 
\cite{bib:krichever-etal94} of the elliptic Calogero-Moser 
system.  Actually, a multi-spin generalization of 
this construction is also possible, which is nothing 
but the genus-one case of Levin and Olshanetsky's 
framework.  

This paper is organized as follows.  In Section 2, 
we illustrate our construction of isomonodromic 
systems in the case of the most classical $A_{\ell-1}$ 
model. This will serve as a prototype of the subsequent 
discussion.   Section 3 is devoted to the models treated 
by the root type Lax pairs, and Section 4 to the spin 
generalization.  Section 5 is for concluding remarks.  
Technically complicated calculations are collected in 
Appendices.

%%%%%%%%%%%%%%%%%%%%%%%%%%%%%%%%%%%%%%%%%%%%%%%%%%%%%%%%%%%%%%%%%%%%
\section{Isomonodromic Systems on the Torus --- a Prototype}
\setcounter{equation}{0}

We start with illustrating our construction for the most 
fundamental case --- the the $A_{\ell-1}$ model and its 
Lax pair in the vector representation of $SU(\ell)$. 

\subsection{$A_{\ell-1}$ Model of Elliptic Calogero-Moser Systems} 

The $A_{\ell-1}$ model is defined by the Hamiltonian 
\beqn
    \calH = \frac{1}{2} \sum_{j=1}^\ell p_j^2 
    + \frac{g^2}{2} \sum_{j \not= k} \wp(q_j - q_k). 
\eeqn
Here $q_j$ and $p_j$ ($j=1,\cdots,\ell$) are the 
coordinates and momenta of the particles with the 
canonical Poisson brackets 
\beqn
    \{q_j,p_k\} = \delta_{jk}, \quad 
    \{q_j,q_k\} = \{p_j,p_k\} = 0. 
\eeqn
Following Manin's equation, we noralize the primitive 
periods as 
\beqn
    2\omega_1 = 1, \quad 2\omega_3 = \tau
\eeqn
The equations of motion are 
give by the canonical equations 
\beqn
    \frac{dq_j}{dt} &=& \{q_j, \calH\} 
      = p_j, 
    \nonumber \\
    \frac{dp_j}{dt} &=& \{p_j, \calH\} 
      = - g^2 \sum_{k \not= j} \wp'(q_j - q_k). 
\eeqn

This elliptic Calogero-Moser system has a Lax pair 
of the form 
\beqn
    L(z) &=& \sum_{j=1}^\ell p_j E_{jj} 
      + ig \sum_{j \not= k} x(q_j - q_k, z) E_{jk},
    \nonumber \\
    M(z) &=& \sum_{j=1}^\ell D_j E_{jj} 
      + ig \sum_{j \not= k} y(q_j - q_k, z) E_{jk}, 
\eeqn
where $E_{jk}$ is the matrix unit, 
$(E_{jk})_{mn} = \delta_{mj}\delta_{nk}$. 
The diagonal elements $D_j$ of $M(z)$ are given by 
\beqn
    D_j = ig \sum_{k \not= j} \wp(q_j - q_k), 
\eeqn
and $x(u,z)$ is a function that satisfies, along 
with its $u$-derivative 
\beqn
    y(u,z) = \frac{\rd x(u,z)}{\rd u}, 
\eeqn
the functional equations 
\beqn
    && x(u,z) y(v,z) - y(u,z) x(v,z) 
       = x(u+v,z) \bigl( \wp(u) - \wp(v) \bigr), 
         \label{eq:sum-rule} \\
    && x(u,z) y(-u,z) - y(u,z) x(-u,z) 
       = \wp'(u), 
         \label{eq:zero-sum-rule} \\
    && x(u,z) x(-u,z) = \wp(z) - \wp(u). 
         \label{eq:factor}
\eeqn
Using these functional equations, one can easily prove 
the following well known result \cite{bib:krichever80}: 
\begin{proposition} The matrices $L(z)$ and $M(z)$ 
satisfy the Lax equation 
\beqn
    \frac{\rd L(u)}{\rd t} = [L(z), M(z)]. 
\eeqn
\end{proposition}
As far as the elliptic Calogero-Moser system is 
concerned, the choice of $x(u,z)$ and $y(u,y)$ 
is rather irrelevant.  A standard choice  is 
the function 
\beqn
    x(u,z) = \frac{\sigma(z - u)}{\sigma(z) \sigma(u)}, 
    \label{eq:x-sigma}
\eeqn
where $\sigma(u) = \sigma(u \mid 1,\tau)$ is the 
Weierstrass sigma function with primitive periods 
$1$ and $\tau$. 

Thus, the elliptic Calogero-Moser system is an 
isospectral integrable system.  An involutive 
set of conserved quantities can be extracted 
from the traces $\Tr L(z)^k$, $k = 2,3,\cdots$ 
of powers of the L-matrix.  The quadratic trace 
is substantially the Hamiltonian 
itself:
\beqn
    \Tr \frac{L(z)^2}{2} 
    = \calH + (\mbox{independent of $p$ and $q$}). 
\eeqn

The functions $x(u,z)$ and $y(u,z)$ based on the 
sigma function, however, are not very suited for 
constructing an isomonodromic system.  We shall 
show an alternative in the next subsection.

\subsection{Our choice of $x(u,z)$ and $y(u,z)$}

Inspired by the work of Levin and Olshanetsky 
\cite{bib:levin-olsh97}, we take the following 
function $x(u,z)$ and its $u$-derivative $y(u,z)$ 
for constructing an isomonodromic Lax pair: 
\beqn
    x(u,z) = \frac{\theta_1(z - u) \theta_1'(0)}
                  {\theta_1(z) \theta_1(u)}. 
    \label{eq:x-theta} 
\eeqn
Here $\theta_1(u)$ is one of Jacobi's elliptic theta 
functions, 
\beqn
    \theta_1(u) = \theta_1(u \mid \tau) 
    = - \sum_{n=-\infty}^{\infty} \exp\left( 
      \pi i\tau \left(n + \frac{1}{2}\right)^2 
      + 2 \pi i \left(n + \frac{1}{2}\right) 
                \left(u + \frac{1}{2}\right) \right), 
\eeqn
and $\theta_1'(u)$ its derivative.  Accordingly, 
the partner $y(u,z)$ can be written 
\beqn
    y(u,z) = - x(u,z) \bigl( \rho(u) + \rho(z - u) \bigr), 
\eeqn
where $\rho(u)$ denotes the logarithmic derivative of 
$\theta_1(u)$, 
\beqn
    \rho(u) = \frac{\theta_1'(u)}{\theta_1(u)}. 
\eeqn
The function $\rho(u)$, too, plays an important role 
throughout this paper.  

\begin{proposition}
These functions $x(u,z)$ and $y(u,z)$ satisfy the functional 
equations (\ref{eq:sum-rule}) -- (\ref{eq:factor}) and 
the differential equation 
\beqn
      2 \pi i \frac{\rd x(u,z)}{\rd \tau} 
      + \frac{\rd^2 x(u,z)}{\rd u \rd z} = 0. 
         \label{eq:x-heat-eq}
\eeqn
\end{proposition}

The last differential equation  (a kind of 
$1+2$-dimensional ``heat equation'')  is a 
characteristic of our $(x,y)$ pair, and plays 
a key role in our construction of isomonodromic 
systems.  

We give a proof of these properties in Appendix A. 
The following are supplementary remarks on these 
functions.  

\begin{itemize}
\item 
The proof of (\ref{eq:sum-rule}--\ref{eq:factor}) 
is based on the following analytical properties of 
$x(u,z)$: 
\begin{enumerate}
\item 
$x(u,z)$ is a meromorphic function of $u$ and $z$.  
The poles on the $u$ plane and the $z$ plane are 
both located at the lattice points $u = m + n\tau$ 
and $z = m + n\tau$ ($m,n \in \bbZ$).  
\item
$x(u,z)$ has the following quasi-periodicity: 
\beqn
    x(u + 1, z) = x(u,z), && 
    x(u + \tau, z) = e^{2\pi i z} x(u,z), 
    \nonumber\\
    x(u, z + 1) = x(u,z), &&
    x(u, z + \tau) = e^{2\pi i u} x(u,z). 
\eeqn
\item
At the origin of the $u$ and $z$ planes, 
$x(u,z)$ exhibits the following singular behavior: 
\beqn
    x(u,z) &=& \frac{1}{u} - \rho(z) + O(u) 
      \quad (u \to 0), 
    \nonumber \\
    x(u,z) &=& - \frac{1}{z} + \rho(u) + O(z) 
      \quad (z \to 0). 
\eeqn
\end{enumerate}

\item 
These properties are an immediate consequence 
of the following well known fact: 

\begin{enumerate}
\item 
$\theta_1(u)$ is an entire function with simple zeros 
at the lattice points $u = m + n\tau$ ($m,n \in \bbZ$). 
\item 
$\theta_1(u)$ is an odd and quasi-periodic function, 
\beqn
    && \theta_1(-u) = \theta_1(u + 1) 
      = - \theta_1(u), 
    \nonumber \\
    && \theta_1(u + \tau) 
      = - e^{-\pi i \tau - 2 \pi i u} \theta_1(u). 
\eeqn
\end{enumerate}

\item 
One can similarly see the following analytical 
properties of $\rho(u)$: 
\begin{enumerate}
\item
$\rho(u)$ is a meromorphic function with poles at 
the lattice points $u = m + n\tau$ ($m,n \in \bbZ$). 
\item
$\rho(u)$ is an odd function with additive 
quasi-periodicity: 
\beqn
    \rho(-u) = - \rho(u), \quad 
    \rho(u + 1) = \rho(u), \quad 
    \rho(u + \tau) = \rho(u) - 2 \pi i. 
\eeqn
\item
At the origin $u = 0$, $\rho(u)$ exhibits the 
following singular behavior: 
\beqn
    \rho(u) = \frac{1}{u} 
    + \frac{\theta_1'''(0)}{3\theta_1'(0)} u 
    + O(u^3) 
    \quad (u \to 0). 
\eeqn
\end{enumerate}

\item 
The proof of (\ref{eq:x-heat-eq}) is based 
on the well known ``heat equation'' 
\beqn
    4 \pi i \frac{\rd \theta_1(u)}{\rd \tau} 
      = \frac{\rd^2 \theta_1(u)}{\rd u^2}. 
    \label{eq:theta-heat-eq} 
\eeqn
of the Jacobi theta function. 
\end{itemize}

\subsection{Isomonodromic deformations} 

Replacing $d/dt \to 2 \pi id/d\tau$, one obtains 
a non-autonomous Hamiltonian system: 
\beqn
    && 2 \pi i \frac{dq_j}{d\tau} = \{q_j, \calH\} 
       = p_j, 
    \nonumber \\
    && 2 \pi i \frac{dp_j}{d\tau} = \{p_j, \calH\} 
       = - g^2 \sum_{k \not= j} \wp'(q_j - q_k). 
\eeqn
We now demonstrate that this gives an isomonodromic system 
on the torus $E_\tau$.  A key is the following Lax equation: 

\begin{proposition} $L(z)$ and $M(z)$ satisfy 
the Lax equation 
\beqn
    2 \pi i \frac{\rd L(z)}{\rd \tau} 
    + \frac{\rd M(z)}{\rd z} 
    = [L(z), M(z)]. 
\eeqn
\end{proposition}

\proof 
Let us notice that the right hand side of the isospectral 
Lax equation is in fact the Poisson bracket of $L(z)$ and 
the Hamiltonian: 
\beqn 
    [L(z), M(z)] = \frac{\rd L(z)}{\rd t} = \{L(z), \calH\}. 
\eeqn
Since the phase space and the Hamiltonian are the same 
as those of the original system, the relation 
$[L(z),M(z)] = \{L(z), \calH\}$ persists in the 
present setup.  Thus the right hand side of the 
Lax equation can be written 
\beqn
    [L(z), M(z)] &=& \{L(z), \calH\} 
    \nonumber \\
    &=&  \sum_{j=1}^\ell \{p_j, \calH\} E_{jj} 
    + ig \sum_{j \not= k} \{q_j - q_k, \calH\} 
          y(q_j - q_k, z) E_{jk}. 
\eeqn
On the other hand, 
\beqn
    2 \pi i \frac{\rd L(z)}{\rd \tau} 
    + \frac{\rd M(z)}{\rd z} 
    &=& \sum_{j=1}^\ell 2 \pi i \frac{dp_j}{d\tau} E_{jj} 
    \nonumber \\
    && + ig \sum_{j \not= k} 2 \pi i \left( 
         \frac{dq_j}{d\tau} - \frac{dq_k}{d\tau} \right) 
         y(q_j - q_k, z) E_{jk} 
    \nonumber \\
    && + ig \sum_{j \not= k} 
        \left( 2 \pi i \frac{\rd x(u,z)}{\rd \tau} 
        + \frac{\rd y(u,z)}{\rd z} \right)_{u = q_j - q_k} E_{jk}.     
\eeqn
The last sum vanishes because of the ``heat equation'' 
(\ref{eq:x-heat-eq}).  The other part coincides, 
term-by-term, with the above expression of the commutator 
$[L(z),M(z)]$. 
\qed

This Lax equation enables us to interpret the 
non-autonomous Hamiltonian system as an isomonodromic 
system on the torus $E_\tau$.  The Lax equation is 
nothing but the Frobenius integrability condition of 
a linear system of the form 
\beqn
    \frac{\rd Y(z)}{\rd z} = L(z) Y(z), \quad 
    2 \pi i \frac{\rd Y(z)}{\rd \tau} + M(z) Y(z) = 0. 
\eeqn
The first equation is an ordinary differential equation 
on the torus $E_\tau$, and has a regular singular point 
at $z = 0$.  Analytic continuation of the solution around 
this singular point yields a monodromy matrix $\Gamma_0$.  
Besides this local monodromy matrix, there are global 
monodromy matrices $\Gamma_\alpha$ and $\Gamma_\beta$ 
that arise in analytic continuation along the $\alpha$ 
($z \to z + 1$) and $\beta$ ($z \to z + \tau$) cycles.  
The second equation of the above linear system implies 
that these monodromy matrices are left invariant as 
$\tau$ varies.  

Let us specify this observation in more detail.  
The situation is more complicated than isomonodromic 
systems on the Riemann sphere: The monodromy of 
$L(z)$ and $M(z)$ themselves are non-trivial, 
\beqn
    && L(z + 1) = L(z), \quad M(z + 1) = M(z), 
    \nonumber \\
    && L(z + \tau) = e^{2\pi iQ} L(z) e^{- 2\pi iQ}, 
    \nonumber \\ 
    && M(z + \tau) = e^{2\pi iQ} \bigl( 
       M(z) + 2\pi i L(z) \bigr) e^{-2\pi iQ} 
       - 2\pi iP, 
\eeqn
where $Q = \sum_{j=1}^\ell q_j E_{jj}$ and 
$P = \sum_{j=1}^\ell p_j E_{jj}$.  These relations 
are a consequence of the quasi-periodicity of 
$x(u,z)$, $y(u,z)$ and $\rho(z)$.  
The monodromy of $L(z)$ implies that $Y(z)$ 
has to be treated as a section of a non-trivial 
$GL(\ell,\bbC)$-bundle (or $SL(\ell,\bbC)$-bundle, 
if we take the center of mass frame with 
$\sum_{j=1}^\ell p_j = 0$) on the torus $E_\tau$.  
The monodromy matrices $\Gamma_0$, $\Gamma_\alpha$ 
and $\Gamma_\beta$ thus arise as follows: 
\beqn
    Y(z e^{2\pi i}) = Y(z) \Gamma_0, \quad 
    Y(z + 1) = Y(z) \Gamma_\alpha, \quad 
    Y(z + \tau) = e^{2\pi iQ} Y(z) \Gamma_\beta. 
\eeqn
Note that the exponential factor in the last relation 
reflects the non-trivial monodromy of $L(z)$ along 
the $\beta$-cycle.  
Having this monodromy structure of $Y(z)$, 
one can deduce the following fundamental observation: 

\begin{proposition} The monodromy matrices 
do not depend on $\tau$, i.e., 
\beqn
    \frac{d\Gamma_0}{d\tau} 
    = \frac{d\Gamma_\alpha}{d\tau}
    = \frac{d\Gamma_\beta}{d\tau} 
    = 0. 
\eeqn
\end{proposition}

\proof
Let us rewrite the second equation of the linear system as 
\beqn
    M(z) = - 2\pi i \frac{\rd Y(z)}{\rd \tau} Y(z)^{-1}, 
\eeqn
and examine the implication of the monodromy structure 
of $Y(z)$ noted above.  This leads to the following relations: 
\beqn
    M(z e^{2\pi i}) &=& M(z) 
      - 2\pi i Y(z) \frac{\rd\Gamma_0}{\rd\tau} 
        \Gamma_0^{-1} Y(z)^{-1}, 
    \nonumber \\
    M(z + 1) &=& M(z) 
      - 2\pi i Y(z) \frac{\rd\Gamma_\alpha}{\rd\tau} 
        \Gamma_\alpha^{-1} Y(z)^{-1}, 
    \nonumber \\
    M(z + \tau) &=& 
      e^{2\pi iQ} \bigl( M(z) + 2 \pi i L(z) \bigr) e^{-2\pi iQ} 
      - 2\pi i P 
    \nonumber \\
    && - 2\pi i Y(z) \frac{\rd\Gamma_\beta}{\rd\tau} 
         \Gamma_\beta^{-1} Y(z)^{-1}. 
\eeqn
(We have used the relation $2\pi i dQ/d\tau = P$.)  
These relations are consistent with the aforementioned 
monodromy structure of $M(z)$ if and only if the monodromy 
matrices of $Y(z)$ are independent of $\tau$. \qed

%%%%%%%%%%%%%%%%%%%%%%%%%%%%%%%%%%%%%%%%%%%%%%%%%%%%%%%%%%%%
\section{Elliptic Calogero-Moser Systems Based on Root Systems}
\setcounter{equation}{0}

Here we consider the elliptic Calogero-Moser systems 
associated with a general irreducible (but not necessary 
reduced) root system $\Delta$.  

In the following, the root system $\Delta$ is assumed to 
be realized in an $\ell$-dimensional Euclidean space 
$M = \bbR^\ell$.  
Let $x \cdot y$ denote the inner product of two vectors 
in $M$ and its bilinear extension to the complexification 
$M^\bbC = M \otimes_\bbR \bbC$. The dual space 
$M^* = Hom(M,\bbR)$ of $M$ is identified with $M$ by 
this inner product.  Each element $\alpha \in \Delta$ 
induces a reflection (the Weyl reflection) $s_\alpha(x) = 
x - (2 \alpha \cdot x / \alpha\cdot\alpha) \alpha$. 
This gives a representation of the Weyl group $W(\Delta)$ 
on $M$.  The root system $\Delta$ is invariant under the 
action of this Weyl group.  

The elliptic Calogero-Moser system associated with 
the root system $\Delta$ is a Hamiltonian system on 
$M \times M$ (or its complexification $M^\bbC \times M^\bbC$).  
The orthognal coordinates 
$(q,p) = (q_1,\cdots,q_\ell,p_1,\cdots,p_\ell)$ of 
$M \times M$ give canonical coordinates and momenta with 
the Poisson brackets 
\beqn
    \{q_j,p_k\} = \delta_{jk}, \quad 
    \{p_j,p_k\} = \{q_j,q_k\} = 0. 
\eeqn

\subsection{Simply laced models}

We first consider the case of simply laced 
($A_{\ell-1}$, $D_\ell$ and $E_\ell$) root systems. 
The associated elliptic Calogero-Moser 
system is defined by the Hamiltonian 
\beqn
    \calH = \frac{1}{2} p \cdot p 
            + \frac{g^2}{2}\sum_{\alpha \in \Delta} 
              \wp(\alpha \cdot q). 
\eeqn
Here $g$ is a coupling constant, and $\wp(u)$ the 
Weierstrass $\wp$ function with primitive periods 
$1$ and $\tau$.  The equations of motion can be 
written 
\beqn
    \frac{dq}{dt} = p, \quad 
    \frac{dp}{dt} = - \frac{g^2}{2} \sum_{\alpha\in\Delta} 
      \wp'(\alpha\cdot q) \alpha. 
\eeqn

We first review the ``root type'' Lax pair of Bordner 
et al. for these models \cite{bib:bordner-etal98a}, 
then explain how to convert these isospectral systems 
to isomonodromic systems.

\subsubsection{Root type Lax pair} 
The ``root type'' Lax pair for these simply laced models 
are $\Delta \times \Delta$ matrices, i.e., matrices 
whose rows and columns are indexed by the root system 
$\Delta$.  They are made of three parts: 
\beqn
    L(z) = P + X_1(z) +X_2(z), \quad 
    M(z) = D + Y_1(z) + Y_2(z). 
\eeqn
$P$ and $D$ are diagonal matrices, 
\beqn
    P_{\beta\gamma} = p \cdot \beta \delta_{\beta\gamma}, 
    \quad 
    D_{\beta\gamma} = D_\beta \delta_{\beta\gamma} 
    \quad (\beta,\gamma \in \Delta), 
\eeqn
and the diagonal elements $D_\beta$ of $D$ are given by 
\beqn
    D_\beta = ig \wp(\beta \cdot q) 
      + ig \sum_{\gamma\in\Delta,\beta\cdot\gamma=1} 
        \wp(\gamma \cdot q). 
\eeqn
$X_1(z)$, etc. are diagonal-free matrices of the form 
\beqn
    X_1(z) &=& ig \sum_{\alpha\in\Delta} 
      x(\alpha \cdot q, z) E(\alpha), 
    \nonumber \\
    X_2(z) &=& 2ig \sum_{\alpha\in\Delta} 
      x(\alpha \cdot q, 2z) E(2\alpha), 
    \nonumber \\
    Y_1(z) &=& ig \sum_{\alpha\in\Delta} 
      y(\alpha \cdot q, z) E(\alpha), 
    \nonumber \\
    Y_2(z) &=& ig \sum_{\alpha\in\Delta} 
      y(\alpha \cdot q, 2z) E(2\alpha), 
\eeqn
where $x(u,z)$ and $y(u,z)$ are the same as the functions 
used in the previous section, and $E(\alpha)$ and 
$E(2\alpha)$ are $\Delta\times\Delta$ matrices of the form 
\beqn
    E(\alpha)_{\beta\gamma} = \delta_{\alpha,\beta-\gamma}, 
    \quad 
    E(2\alpha)_{\beta\gamma} = \delta_{2\alpha,\beta-\gamma}
    \quad (\beta, \gamma \in \Delta). 
\eeqn 
(We have slightly modified the notation of Bordner et al: 
$x(u,2z)$, $y(u,2z)$ and $E(2\alpha)$ amount to 
$x_d(u,z)$, $y_d(u,z)$ and $E_d(\alpha)$ in their notation.) 

These matrices satisfy the Lax equation
\beqn
    \frac{\rd L(z)}{\rd t} = [L(z), M(z)] 
\eeqn
under the equations of motions.  The traces $\Tr L(z)^k$, 
$k = 2,3,\cdots$, of powers of $L(z)$ are conserved, 
and an involutive set of conserved quantities can be 
extracted from these traces.  The Hamiltonian itself 
can be reproduced from the quadratic trace $\Tr L(z)^2$. 
We refer the details of these results to the paper of 
Bordner et al. \cite{bib:bordner-etal98a}.  The choice 
of $x(u,z)$ and $y(u,z)$ is irrelevant in this case, too.  

Thus, in particular, the $A_{\ell-1}$ model turns out 
to have at least two distinct Lax pairs --- the Lax pair 
of $\ell\times\ell$ matrices realized in the vector 
representation of $sl(\ell)$, and the Lax pair of 
$\ell(\ell-1) \times \ell(\ell-1)$ matrices based 
on the $A_{\ell-1}$ root system.  This is also the 
case for the other simply laced root systems. 
Bordner et al. call the Lax pairs of the first type 
the ``minimal type'', because they are realized in 
a minimal representation of the associated (not 
necessary simply laced) Lie algebra.  It should 
be noted that the ``root type'' Lax pairs do not 
possess a Lie algebraic structure; unlike the usual 
root basis of simple Lie algebras, the matrices 
$E(\alpha)$ and $E(2\alpha)$ are not closed under 
the Lie bracket.

\subsubsection{Isomonodromic system} 
The prescription for constructing an isomonodromic 
analogue is the same as the previous case, namely, 
to replace $d/dt \to 2 \pi id/d\tau$.  This 
converts the equations of motion of the elliptic 
Calogero-Moser system to the non-autonomous system 
\beqn
    2 \pi i \frac{dq}{dt} = p, \quad 
    2 \pi i \frac{dp}{dt} = 
      - \frac{g^2}{2} \sum_{\alpha\in\Delta} 
      \wp'(\alpha\cdot q) \alpha. 
\eeqn
Let $x(u,z)$ be the function defined in (\ref{eq:x-theta}), 
and $y(u,z)$ its $u$-derivative.  The following 
are the keys to an isomonodromic interpretation. 

\begin{proposition} 
1. $L(z)$ and $K(z)$ satisfy the Lax equation 
\beqn
    2 \pi i \frac{\rd L(z)}{\rd \tau}
      + \frac{\rd M(z)}{\rd z} 
       = [L(z), M(z)] . 
\eeqn
2. $L(z)$ and $M(z)$ have the following monodromy 
property: 
\beqn
    && L(z + 1) = L(z), \quad M(z + 1) = M(z), 
    \nonumber \\
    && L(z + \tau) = e^{2\pi iQ} L(z) e^{-2\pi iQ}, 
    \nonumber \\
    && M(z + \tau) = e^{2\pi iQ} \bigl( M(z) 
       + 2 \pi i L(z) \bigr) e^{-2\pi iQ} - 2 \pi i P, 
\eeqn
where $Q$ is the diagonal matrix with matrix elements 
$Q_{\beta\gamma} = q\cdot \beta \delta_{\beta\gamma}$.  
\end{proposition}

\proof 
The proof is almost the same as the proof for the 
isomonodromic Lax pair of the $A_{\ell-1}$ model in the 
vector representation. Let us first verify the Lax equation. 
The right hand side of the Lax equation can be written 
\beqn
    [L(z), M(z)] 
    &=&  \{P, \calH\} 
       + ig \sum_{\alpha\in\Delta} \{\alpha\cdot q, \calH\} 
          y(\alpha\cdot q, z) E(\alpha) 
    \nonumber \\
    && + 2ig \sum_{\alpha\in\Delta} \{\alpha\cdot q, \calH\} 
          y(\alpha\cdot q, 2z) E(2\alpha). 
\eeqn
On the other hand, 
\beqn
    2\pi i\frac{\rd L(z)}{\rd \tau} + \frac{\rd M(z)}{\rd z} 
    &=&  2 \pi i \frac{\rd P}{\rd \tau} 
       + ig \sum_{\alpha\in\Delta} 
           2\pi i\frac{\rd \alpha\cdot q}{\rd \tau} 
           y(\alpha\cdot q, z) E(\alpha) 
    \nonumber \\
    && + 2ig \sum_{\alpha\in\Delta} 
           2\pi i\frac{\rd \alpha\cdot q}{\rd \tau} 
           y(\alpha\cdot q, 2z) E(2\alpha) 
    \nonumber \\
    && + ig \sum_{\alpha\in\Delta}
           \left( 2\pi i\frac{\rd x(u,z)}{\rd \tau} 
           + \frac{\rd y(u,z)}{\rd z} 
           \right)_{u=\alpha\cdot q} E(\alpha) 
    \nonumber \\
    && + 2ig \sum_{\alpha\in\Delta} 
           \left( 4\pi i\frac{\rd x(u,2z)}{\rd \tau} 
           + \frac{\rd y(u,2z)}{\rd z} 
           \right)_{u=\alpha\cdot q} E(2\alpha). 
\eeqn
The last two sums vanish because of (\ref{eq:x-heat-eq}). 
The other part coincides by the equations of motion. 
Thus we obtain the Lax equation.  Let us next consider 
the monodromy of $L(z)$ and $M(z)$. Note the commutation 
relations 
\beqn
    [Q, E(\alpha)] = q \cdot \alpha E(\alpha), \quad 
    [Q, E(2\alpha)] = 2q \cdot \alpha E(2\alpha), 
\eeqn
which can be exponentiated as follows: 
\beqn
    e^{2\pi iQ} E(\alpha) e^{-2\pi iQ} 
      = e^{2\pi iq\cdot\alpha} E(\alpha), \quad 
    e^{2\pi iQ} E(2\alpha) e^{-2\pi iQ} 
      = e^{4\pi iq\cdot\alpha} E(2\alpha). 
\eeqn
The monodromy property of $L(z)$ and $M(z)$ 
can be derived from these relations and the 
quasi-periodicity of $x(u,z)$ and $y(u,z)$.  \qed

The rest is parallel to the case in the previous 
section.  The only difference is that the ordinary 
differential equation 
\beqn
    \frac{dY(z)}{dz} = L(z) Y(z) 
\eeqn
on the torus $E_\tau$ has {\it four} regular singular 
points at $z = 0,\omega_1,\omega_2,\omega_3$.  The 
latter three singular points originates in $X_2(z)$.  
Let $\Gamma_a$ ($a = 0,1,2,3$) denote the monodromy 
matrices in analytic continuation of $Y(z)$ around 
these four points.  The Lax equation implies that 
these local monodromy matrices and the two global 
ones $\Gamma_\alpha$ and $\Gamma_\beta$ are 
independent of $\tau$: 
\beqn
      \frac{\rd \Gamma_0}{\rd \tau} 
    = \cdots 
    = \frac{\rd \Gamma_3}{\rd \tau} 
    = \frac{\rd \Gamma_\alpha}{\rd \tau}
    = \frac{\rd \Gamma_\beta}{\rd \tau} 
    = 0. 
\eeqn

\subsection{Non-simply laced models}

The elliptic Calogero-Moser system associated with 
a non-simply laced ($B_\ell$, $C_\ell$, $F_4$, $G_2$ 
and $BC_\ell$) root systems can have several independent 
coupling constants, one for each Weyl group orbit in 
the root system.  The root type Lax pairs are extended 
to the non-simply laced cases by Bordner et al. 
\cite{bib:bordner-etal98b}. As they pointed out, 
one can construct a different root type Lax pair 
for each Weyl group orbit of the root system. Thus 
the $B_\ell$, $C_\ell$, $F_4$ and $G_2$ models have, 
respectively, two distinct Lax pairs based on the 
orbits of long and short roots, whereas the $BC_\ell$ 
model has three based on the orbits of long, middle, 
and short roots.  Note that each Weyl group orbit 
consists of roots of the same length. 

Although all the non-simply laced models can be treated 
in the same way, let us illustrate our construction of 
isomonodromic systems for the $BC_\ell$ model.  This is 
also intended to be a prototype of the case that we 
shall consider in the next subsection.  

\subsubsection{$BC_\ell$ model} 
The $BC_\ell$ root system can be realized in 
$M = \bbR^\ell$: 
\beqn
    \Delta(BC_\ell) &=& 
      \Delta_l \cup \Delta_m \cup \Delta_s, 
    \nonumber \\
    \Delta_l &=& \{ \pm 2e_j \mid 1 \le j \le \ell \} 
      \quad (\mbox{long roots}), 
    \nonumber \\
    \Delta_m &=& \{ \pm e_j \pm e_k \mid j \not= k \}
      \quad (\mbox{middle roots}), 
    \nonumber \\
    \Delta_s &=& \{ \pm e_j \mid 1 \le j \le \ell \} 
      \quad (\mbox{short roots}), 
\eeqn
where $e_1,\cdots,e_\ell$ are the standard orthonormal 
basis of $\bbR^\ell$.  $\Delta_l$, $\Delta_m$ and 
$\Delta_s$ give the three Weyl group orbits. 

The Hamiltonian of the $BC_\ell$ model takes the form 
\beqn
    \calH = \frac{1}{2} p \cdot p 
    + \frac{g_m^2}{2} \sum_{\alpha\in\Delta_m} 
        \wp(\alpha \cdot q) 
    + \frac{g_l^2}{4} \sum_{\alpha\in\Delta_l}
        \wp(\alpha \cdot q) 
    + \gtilde_s^2 \sum_{\alpha\in\Delta_s} 
        \wp(\alpha \cdot q). 
\eeqn
The equations of motion can be written 
\beqn
    \frac{dq}{d\tau} &=& p, 
    \nonumber \\
    \frac{dp}{d\tau} &=& 
      - \frac{g_m^2}{2} \sum_{\alpha\in\Delta_m} 
          \wp'(\alpha\cdot q) \alpha 
      - \frac{g_l^2}{4} \sum_{\alpha\in\Delta_l} 
          \wp'(\alpha\cdot q) \alpha 
      - \gtilde_s^2 \sum_{\alpha\in\Delta_s} 
          \wp'(\alpha\cdot q) \alpha. 
\eeqn
$g_m,g_l$ and $\gtilde_s$ are three independent coupling 
constants.  $\gtilde_s$ is a modified (``renormalized'' 
in the terminology of Bordner et al.) coupling constant 
connected with a more fundamental (``bare'', so to speak) 
coupling constant $g_s$ as 
\beqn 
    \gtilde_s^2 = g_s^2 + \frac{g_s g_l}{2}. 
\eeqn
The ``bare'' coupling constant appears in the 
construction of a Lax pair.

\subsubsection{Root type Lax pair for $BC_\ell$ model} 
As mentioned above, there are at least three root type 
Lax pairs based on the three Weyl group orbits 
$\Delta_m$, $\Delta_l$ and $\Delta_s$.  Bordner et al. 
constructed only one of them,  namely, a Lax pair based 
on $\Delta_m$. Here we present a Lax pair based on 
$\Delta_s$.  This is a $2\ell \times 2\ell$ system, 
much smaller than the Lax pair based on $\Delta_m$, 
and presumably more suitable for studying the associated 
isomonodromic deformations. 

The Lax pair are indexed by $\Delta_s$ and take the 
following form: 
\beqn
    L(z) &=& P + X_1(z) + X_2(z) + X_3(z), 
    \nonumber \\
    M(z) &=& D + Y_1(z) + Y_2(z) + Y_3(z). 
\eeqn
$P$ and $D$ are diagonal matrices, 
\beqn
    P_{\beta\gamma} = p \cdot \beta \delta_{\beta\gamma}, 
    \quad 
    D_{\beta\gamma} = D_\beta \delta_{\beta\gamma} 
    \quad (\beta,\gamma \in \Delta_s), 
\eeqn
and the diagonal elements of $D$ are given by 
\beqn
    D_\beta 
    = ig_m \sum_{\gamma\in\Delta_s,\beta\cdot\gamma=1} 
      \wp(\gamma\cdot q) 
    + ig_l \wp(2\beta\cdot q) + ig_s \wp(\beta\cdot q). 
\eeqn
$X_1(z)$, etc. are diagonal-free matrices of the form 
\beqn
    X_1(z) &=& ig_m \sum_{\alpha\in\Delta_m} 
      x(\alpha\cdot q, z) E(\alpha), 
    \nonumber \\
    X_2(z) &=& ig_l \sum_{\alpha\in\Delta_l} 
      x(\alpha\cdot q, z) E(\alpha), 
    \nonumber \\
    X_3(z) &=& 2ig_s \sum_{\alpha\in\Delta_s}  
      x(\alpha\cdot q, 2z) E(2\alpha), 
    \nonumber \\
    Y_1(z) &=& ig_m \sum_{\alpha\in\Delta_m} 
      y(\alpha\cdot q, z) E(\alpha), 
    \nonumber \\
    Y_2(z) &=& ig_l \sum_{\alpha\in\Delta_l} 
      y(\alpha\cdot q, z) E(\alpha), 
    \nonumber \\
    Y_3(z) &=& ig_s \sum_{\alpha\in\Delta_s} 
      y(\alpha\cdot q, 2z) E(2\alpha), 
\eeqn
where 
\beqn
    E(\alpha)_{\beta\gamma} = \delta_{\alpha,\beta-\gamma}, 
    \quad 
    E(2\alpha)_{\beta\gamma} = \delta_{2\alpha,\beta-\gamma} 
    \quad (\beta,\gamma \in \Delta_s). 
\eeqn
This Lax pair is a specialization of the Lax pair for 
the extended twisted model that we shall present in 
the next subsection.

\subsubsection{Isomonodromic system} 
This system, too, can be converted to an isomonodromic 
system by replacing $d/dt \to 2\pi id/d\tau$. The equations 
of motion are a non-autonomous system of the form 
\beqn
    2 \pi i \frac{dq}{d\tau} &=& p, 
    \nonumber \\
    2 \pi i \frac{dp}{d\tau} &=& 
      - \frac{g_m^2}{2} \sum_{\alpha\in\Delta_m} 
          \wp'(\alpha\cdot q) \alpha 
      - \frac{g_l^2}{4} \sum_{\alpha\in\Delta_l} 
          \wp'(\alpha\cdot q) \alpha 
      - \gtilde_s^2 \sum_{\alpha\in\Delta_s} 
          \wp'(\alpha\cdot q) \alpha. 
\eeqn
The following can be verified just as in the case of 
simply lased models: 
\begin{enumerate}
\item 
$L(z)$ and $M(z)$ satisfy the Lax equation 
\beqn
    2 \pi i \frac{\rd L(z)}{\rd \tau} 
    + \frac{\rd M(z)}{\rd z} 
    = [L(z), M(z)]. 
\eeqn 
\item 
$L(z)$ and $M(z)$ have the following monodromy 
property: 
\beqn
    && L(z + 1) = L(z), \quad M(z + 1) = M(z), 
    \nonumber \\
    && L(z + \tau) = e^{2\pi iQ} L(z) e^{-2\pi iQ}, 
    \nonumber \\
    && M(z + \tau) = e^{2\pi iQ} \bigl( M(z) 
       + 2 \pi i L(z) \bigr) e^{-2\pi iQ} - 2 \pi i P. 
\eeqn
\end{enumerate} 
The interpretation of this Lax equation, too, is 
parallel to the simply laced models.  The ordinary 
differential equation 
\beqn
    \frac{dY(z)}{dz} = L(z) Y(z)
\eeqn
on the torus $E_\tau$ has four regular singular points 
at $z = 0,\omega_1,\omega_2,\omega_3$.  The local 
monodromy matrices $\Gamma_a$ ($a = 0,1,2,3$) at these 
points and the global monodromy matrices $\Gamma_\alpha$ 
and $\Gamma_\beta$ are invariant as $\tau$ varies.

\subsection{Twisted and extended twisted models} 

We now proceed to the ``twisted'' and ``extended 
twisted'' models.  The Hamiltonian of the untwisted 
models can be generally written 
\beqn
    \calH = \frac{1}{2} p \cdot p 
      + \frac{1}{2} \sum_{\alpha\in\Delta} 
        g_{|\alpha|}^2 \wp(\alpha\cdot q). 
\eeqn
The twisted models, introduced by D'Hoker and Phong 
\cite{bib:dhoker-phong98} for non-simply laced root 
systems, are defined by a Hamiltonian of the form 
\beqn
    \calH = \frac{1}{2} p \cdot p 
      + \frac{1}{2} \sum_{\alpha\in\Delta} 
        g_{|\alpha|}^2 \wp_{\nu(\alpha)}(\alpha\cdot q), 
\eeqn
where $\wp_{\nu(\alpha)}(u)$ are the $\wp$-functions 
with suitably rescaled primitive periods.  D'Hoker 
and Phong proved the integrability of those twisted 
models by constructing a Lax pair in a representation 
of the associated Lie algebra.  Bordner and Sasaki 
\cite{bib:bordner-etal98c} proposed an alternative 
approach based on root systems rather than Lie 
algebras, and pointed out that the twisted model of 
the $B_\ell$, $C_\ell$ and $BC_\ell$ types can be 
further extended. The extended twisted models have 
one (for the $B_\ell$ and $C_\ell$ models) or two 
(for the $BC_\ell$ model) extra types of elliptic 
potentials.  

Our construction of isomonodromic systems can be 
extended to the twisted and extended twisted 
models.  We illustrate this result, just as in 
the previous subsection, for the $BC_\ell$ model.  
As Bordner and Sasaki noted, the extended twisted 
$BC_\ell$ model is made of five different types 
of elliptic potentials, and coincides with the 
Inozemtsev system \cite{bib:inozemtsev89}.

\subsubsection{Extended twisted $BC_\ell$ model} 
The extended twisted $BC_\ell$ model is defined 
by the Hamiltonian 
\beqn
    \calH 
    &=&  \frac{1}{2} p \cdot p 
       + \frac{g_m^2}{2} \sum_{\alpha\in\Delta_m} 
         \wp(\alpha \cdot q) 
       + \frac{g_{l1}^2}{4} \sum_{\alpha\in\Delta_l} 
         \wp(\alpha \cdot q) 
       + \frac{\gtilde_{l2}^2}{4} \sum_{\alpha\in\Delta_l} 
         \wp^{(2)}(\alpha \cdot q) 
    \nonumber \\
    && + \gtilde_{s1}^2 \sum_{\alpha\in\Delta_s} 
         \wp(\alpha \cdot q) 
       + \gtilde_{s2}^2 \sum_{\alpha\in\Delta_s} 
         \wp^{(1/2)}(\alpha \cdot q).  
\eeqn
$\gtilde_{l2}$, $\gtilde_{s1}$ and $\gtilde_{s2}$ 
are ``renormalized'' coupling constants, which are 
related to unrenormalized coupling constants $g_{l2}$, 
$g_{s1}$ and $g_{s2}$ as follows: 
\beqn
    \gtilde_{l2}^2 &=& g_{l2}^2 + 2 g_{l1} g_{l2}, 
    \nonumber \\
    \gtilde_{s1}^2 &=& g_{s1}^2 + 2 g_{s1} g_{s2} 
      + \frac{1}{2} (g_{s1}g_{l1} + g_{s1}g_{l2} 
        + g_{s2}g_{l2}), 
    \nonumber \\
    \gtilde_{s2}^2 &=& g_{s2}^2 + \frac{g_{s2}g_{l1}}{2}. 
\eeqn
$\wp^{(1/2)}$ and $\wp^{(2)}$ are the $\wp$ functions 
with rescaled primitive periods: 
\beqn
    \wp^{(1/2)}(u) = \wp(u \mid \frac{1}{2}, \tau), 
    \quad 
    \wp^{(2)}(u) = \wp(u \mid 2, \tau). 
\eeqn
(This Hamiltonian is slightly different from the 
Hamiltonian of Bordner and Sasaki, though the 
contents are essentially the same.  With this 
modification, this model reduces to the untwisted 
$BC_\ell$ model as $g_{l2} \to 0$ and $g_{s2} \to 0$.)

\subsubsection{Root type Lax pair for extended twisted 
$BC_\ell$ model}
One can construct, like the untwisted model, three different 
root type Lax pairs can be constructed based on the three 
Weyl group orbits $\Delta_m$, $\Delta_l$ and $\Delta_s$.  
The Lax pair based on $\Delta_m$ is presented by Bordner 
and Sasaki.  The Lax pair based on $\Delta_s$ can be 
obtained by modifying the Lax pair for the untwisted 
$BC_\ell$ model as follows.  

The Lax pair $L(z)$ and $M(z)$ are indexed by $\Delta_s$ 
and made of four parts, 
\beqn
    L(z) &=& P + X_1(z) + X_2(z) + X_3(z), 
    \nonumber \\
    M(z) &=& D + Y_1(z) + Y_2(z) + Y_3(z). 
\eeqn
The diagonal matrix $P$ is the same as the $P$ in the 
untwisted model. The diagonal matrices of $D$ are given 
by 
\beqn
    D_\beta 
    &=&  i g_m \sum_{\gamma\in\Delta_m,\beta\cdot\gamma=1}
               \wp(\gamma \cdot q) 
         + ig_{l1} \wp(2 \beta \cdot q) 
         + ig_{l2} \wp^{(2)}(2 \beta \cdot q) 
    \nonumber \\
    &&   + ig_{s1} \wp(\beta \cdot q) 
         + ig_{s2} \wp^{(1/2)}(\beta \cdot q). 
\eeqn
$X_1(z)$ and $Y_1(z)$ are the same as those for the 
untwisted model.  The other matrices take the 
following form: 
\beqn
    X_2(z) 
    &=& \sum_{\alpha\in\Delta_l} \Bigl( 
          ig_{l1} x(\alpha \cdot q, z) 
        + ig_{l2} x^{(2)}(\alpha \cdot q, z) 
        \Bigr) E(\alpha), 
    \nonumber \\
    X_3(z) 
    &=& \sum_{\alpha\in\Delta_s} \Bigl( 
          2ig_{s1} x(\alpha \cdot q, 2z) 
        + 2ig_{s2} x^{(1/2)}(\alpha \cdot q, 2z) 
        \Bigr) E(2\alpha), 
    \nonumber \\
    Y_2(z) 
    &=& \sum_{\alpha\in\Delta_l} \Bigl( 
          ig_{l1} y(\alpha \cdot q, z) 
        + ig_{l2} y^{(2)}(\alpha \cdot q, z) 
        \Bigr) E(\alpha), 
    \nonumber \\
    Y_3(z) 
    &=& \sum_{\alpha\in\Delta_s} \Bigl( 
          ig_{s1} y(\alpha \cdot q, 2z) 
        + ig_{s2} y^{(1/2)}(\alpha \cdot q, 2z) 
        \Bigr) E(2\alpha). 
\eeqn
This Lax pair reduces to the Lax pair of the untwisted model 
if $g_{l2} = 0$ and $g_{s2} = 0$.  

The new objects arising here are the functions 
$x^{(1/2)}(u,z)$, $x^{(2)}(u,z)$ and their $u$-derivatives 
\beqn
    y^{(1/2)}(u,z) = \frac{\rd x^{(1/2)}(u,z)}{\rd u}, 
    \quad 
    y^{(2)}(u,z) = \frac{\rd x^{(2)}(u,z)}{\rd u}. 
\eeqn
For the consistency of the Lax equation 
\beqn
    \frac{\rd L(z)}{\rd t} = [L(z), M(z)], 
\eeqn
these functions have to satisfy several functional 
equations.  D'Hoker and Phong \cite{bib:dhoker-phong98} 
and Bordner and Sasaki \cite{bib:bordner-etal98c} 
use a set of functions based on the Weierstrass 
sigma functions. We use the function $x(u,z) = 
x(u,z \mid \tau)$ defined in (\ref{eq:x-theta}) 
and its modifications 
\beqn
    x^{(1/2)}(u,z) 
    &=& 2 x(2u, z \mid 2\tau) 
    = \dfrac
      {2\theta_1(z - 2u \mid 2\tau) \theta_1'(0 \mid 2\tau)}
      {\theta_1(z \mid 2\tau) \theta_1(2u \mid 2\tau)}, 
    \nonumber \\
    x^{(2)}(u,z) 
    &=& \frac{1}{2} 
        x(\frac{u}{2}, z \mid \frac{\tau}{2}) 
    = \dfrac
      {\theta_1(z - \frac{u}{2} \mid \frac{\tau}{2}) 
        \theta_1'(0 \mid \frac{\tau}{2})}
      {2\theta_1(z \mid \frac{\tau}{2})
        \theta_1(\frac{u}{2} \mid \frac{\tau}{2})}. 
\eeqn
These functions $x^{(1/2)}(u,z)$ and $x^{(2)}(u,z)$, too, 
satisfy $1+2$-dimensional ``heat equations'' of the form 
\beqn
  && 2 \pi i \frac{\rd x^{(1/2)}(u,z)}{\rd \tau} 
     + \frac{\rd^2 x^{(1/2)}(u,z)}{\rd u \rd z} = 0, 
  \nonumber \\
  && 2 \pi i \frac{\rd x^{(2)}(u,z)}{\rd \tau} 
     + \frac{\rd^2 x^{(2)}(u,z)}{\rd u \rd z} = 0. 
\eeqn
The functional identities for these functions and 
the proof of the Lax equation are presented in 
Appendices B and C.

\subsubsection{Isomonodromic system} 
Replacing $d/dt \to 2\pi id/d\tau$, we obtain a 
non-autonomous Hamiltonian system with the same 
Hamiltonian.  The isomonodromic interpretation 
of this non-autonomous system is again based 
on the following two observations: 
\begin{enumerate} 
\item 
$L(z)$ and $M(z)$ satisfy the Lax equation 
\beqn
    2 \pi i \frac{\rd L(z)}{\rd \tau} 
    + \frac{\rd M(z)}{\rd z} 
    = [L(z), M(z)]. 
\eeqn
\item 
The monodromy of $L(z)$ and $M(z)$ is the same as 
the monodromy of the Lax pair for the untwisted model: 
\beqn
    && L(z + 1) = L(z), \quad M(z + 1) = M(z), 
    \nonumber \\
    && L(z + \tau) = e^{2\pi iQ} L(z) e^{-2\pi iQ}, 
    \nonumber \\
    && M(z + \tau) = e^{2\pi iQ} \bigl( M(z) 
       + 2 \pi i L(z) \bigr) e^{-2\pi iQ} - 2 \pi i P. 
\eeqn
\end{enumerate} 
The ordinary differential equation defined on the torus 
$E_\tau$ by the matrix $L(z)$ has four regular singular 
points at $u = 0,\omega_1,\omega_2,\omega_3$.  The Lax 
equation and the monodromy of $L(z)$ and $M(z)$ ensure 
that the local monodromy matrices $\Gamma_a$ ($a = 0,1,2,3$) 
and the global monodromy matrices $\Gamma_\alpha$ and 
$\Gamma_\beta$ are independent of $\tau$.

\subsubsection{Relation to Inozemtsev system} 
The final task is to clarify the relation to the 
Inozemtsev system.  In terms of the orthogonal 
coordinates $q_j = q \cdot e_j$ and $p_j = p \cdot e_j$ 
($j = 1,\cdots,\ell$), the aforementioned Hamiltonian 
can be written 
\beqn
    \calH &=& \frac{1}{2} \sum_{j=1}^\ell p_j^2 
       + \frac{g_m^2}{2} \sum_{\epsilon,\epsilon'=\pm 1}
         \sum_{j\not= k} \wp(\epsilon q_j + \epsilon' q_k) 
       + \frac{g_{l1}^2}{2} \sum_{j=1}^\ell \wp(2q_j) 
    \nonumber \\
    && + \frac{\gtilde_{l2}^2}{2} 
         \sum_{j=1}^\ell \wp^{(2)}(2q_j) 
       + 2\gtilde_{s1}^2 \sum_{j=1}^\ell \wp(q_j) 
       + 2\gtilde_{s2}^2 \sum_{j=1}^\ell \wp^{(1/2)}(q_j). 
\eeqn
One can rewrite this Hamiltonian using the identities 
\beqn
    \wp(2u) &=& \frac{1}{4} \wp(u) 
       + \frac{1}{4} \wp(u + \omega_1) 
       + \frac{1}{4} \wp(u + \omega_2) 
       + \frac{1}{4} \wp(u + \omega_3), 
    \nonumber \\
    \wp^{(1/2)}(u) &=& \wp(u) 
       + \wp(u + \omega_1) - \wp(\omega_1), 
    \nonumber \\
    \wp^{(2)}(2u) &=& \frac{1}{4}\wp(u) 
       + \frac{1}{4} \wp(u + \omega_3) 
       - \frac{1}{4} \wp(\omega_3). 
\eeqn
The outcome is, up to a term $h(\tau)$ depending on 
$\tau$ only, the Inozemtsev Hamiltonian 
\beqn
    \calH &=& \frac{1}{2}\sum_{j=1}^\ell p_j^2 
       + \frac{g_m^2}{2} \sum_{\epsilon,\epsilon'=\pm 1}
         \sum_{j\not= k} \wp(\epsilon q_j + \epsilon' q_k) 
       + \sum_{j=1}^\ell \sum_{a=0}^3 
           g_a^2 \wp(q_j + \omega_a) 
       + h(\tau). 
\eeqn
The coupling constants $g_a$ $(a=0,1,2,3$) are given by 
\beqn
    g_0^2 = \frac{1}{8} (g_{l1}^2 + \gtilde_{l2}^2) 
          + 2(\gtilde_{s1}^2 + \gtilde_{s2}^2), 
    &&
    g_1^2 = \frac{g_{l1}^2}{8} + 2 \gtilde_{s2}^2, 
    \nonumber \\
    g_2^2 = \frac{g_{l1}^2}{8}, 
    && 
    g_3^2 = \frac{1}{8} (g_{l1}^2 + \gtilde_{l2}^2). 
\eeqn

%%%%%%%%%%%%%%%%%%%%%%%%%%%%%%%%%%%%%%%%%%%%%%%%%%%%%%%%%%%%
\section{Spin Generalization of Elliptic Calogero-Moser Systems}
\setcounter{equation}{0}

``Spin generalization'' is a generalization of the elliptic 
Calogero-Moser systems coupled to spin degrees of freedom.  
Such a spin generalization is characterized by a simple 
Lie algebra rather than a root system.  The (classical) 
spin variables take values in the dual space $\frakg^*$, 
or a coadjoint orbit therein, of the Lie algebra $\frakg$.   
We shall first examine the $sl(\ell)$ model as a prototype, 
then proceed to the models based on a general simple Lie 
algebra.

\subsection{Spin generalization for $sl(\ell)$} 

The $sl(\ell)$ spin generalization was first introduced by 
Krichever et al. \cite{bib:krichever-etal94}. They obtained 
the spin generalization, just like the spinless case 
\cite{bib:krichever80}, via the pole dynamics of the 
matrix KP hierarchy.  

\subsubsection{Hamiltonian formalism} 
This model is a constrained Hamiltonian system.  
The Hamiltonian is given by 
\beqn
    \calH = \frac{1}{2} \sum_{j=1}^\ell p_j^2 
          - \frac{1}{2} \sum_{j \not= k} 
              \wp(q_j - q_k) F_{jk} F_{kj}. 
\eeqn
Here $q_j$ and $p_j$ ($j = 1,\cdots,\ell$) are the 
canonical coordinates and momenta of the Calogero-Moser 
particles, and $F_{jk}$ ($j,k = 1,\cdots,\ell$) a set 
of classical $sl(\ell)$ spin variables, whose Poisson 
brackets are determined by the Kostant-Kirillov Poisson 
structure on the dual space of $sl(\ell)$: 
\beqn
    \{ F_{jk}, F_{mn} \} 
    = \delta_{mk} F_{jn} - \delta_{jn} F_{mk}. 
\eeqn
The equations of motion can be written 
\beqn
    \frac{dq_j}{dt} 
    &=& p_j, \quad  
         \frac{dp_j}{dt} = 
         \sum_{k \not= j} \wp'(q_j - q_k) F_{jk} F_{kj}, 
    \nonumber \\
    \frac{dF_{jk}}{dt} 
    &=& - \sum_{m \not= j} \wp(q_j - q_m) F_{jm} 
        + \sum_{m \not= k} \wp(q_m - q_k) F_{mk}
    \nonumber \\
    &&  - \wp(q_j - q_k) (F_{jj} - F_{kk}).
\eeqn
In particular, the diagonal elements $F_{jj}$ of 
the spin variables are conserved quantities: 
$dF_{jj}/dt = 0$.  Although the Hamiltonian 
does not contain the diagonal elements explicitly, 
they do appear in the equations of motion.  We now 
put the constraints 
\beqn
    F_{jj} = 0 \quad (j = 1,\cdots,\ell). 
    \label{eq:spin-constraint}
\eeqn
These constraints ensure the integrability. 
(Actually, the integrability is retained if the 
constraints are replaced by $F_{jj} = c$,  
$j = 1,\cdots,\ell$, where $c$ is a constant.)

\subsubsection{Lax pair in vector representation} 
The Lax pair of the spinless $A_{\ell-1}$ model 
in the vector representation of $sl(\ell)$ can 
be readily extended to the spin generalization 
as follows: 
\beqn
    L(z) &=& \sum_{j=1}^\ell p_j E_{jj} 
      + \sum_{j \not= k} \sigma(q_j - q_k, z) F_{kj} E_{jk}, 
    \nonumber \\
    M(z) &=& - \sum_{j \not= k} \sigma(q_j - q_k, z) 
      \bigl( \rho(q_j - q_k) + \rho(z - q_j + w_k)\bigr) 
      F_{kj} E_{jk}, 
\eeqn
where 
\beqn
    \rho(u) = \frac{\theta_1'(u)}{\theta_1(u)}, 
    \quad 
    \sigma(u,z) = \frac{\theta_1(u - z) \theta_1'(0)} 
                       {\theta_1(z) \theta_1(u)}. 
\eeqn
It is these functions that Felder and Wieczerkowski 
used in the KZB equation \cite{bib:feld-wiec94}.  
The function $\rho(u)$ is already familiar to us. 
The function $\sigma(u,z)$ is also just a disguise of 
the function $x(u,z)$ that we have used in the preceding 
sections: 
\beqn
    \sigma(u,z) &=& - x(u,z). 
\eeqn
We however dare to retain the notation of Felder and 
Wieczerkowski so as to stress the similarity with their 
work.  In these notations, the aforementioned functional 
identities of $x(u,z)$ and $y(u,z)$ can be rewritten 
\beqn
    && \sigma(u,z) \sigma(v,z) \bigl( \rho(v) 
       + \rho(z - v) - \rho(u) - \rho(z - u) \bigr) 
     = \sigma(u + v, z) \bigl( \wp(u) - \wp(v) \bigr), 
    \label{eq:sigma-sum-rule} \\
    && 2 \sigma(u,z) \sigma(-u,z) \bigl( \rho(u) 
       + \rho(z - u) \bigr) 
     = - \wp'(u), 
    \label{eq:sigma-zero-sum-rule} \\
    && \sigma(u,z) \sigma(-u,z) 
     = \wp(z) - \wp(u). 
    \label{eq:sigma-factor} 
\eeqn       

Using these functional identities, one can derive 
the Lax equation 
\beqn
    \frac{\rd L(z)}{\rd t} = [L(z), M(z)]. 
\eeqn
Note that the constraints (\ref{eq:spin-constraint}) 
are always assumed when we consider the Lax equation. 
Thus the spin generalization, too, is an isospectral 
integrable system. An involutive set of conserved 
quantities obtained from the traces $\Tr L(z)^k$, 
$k=2,3,,\cdots$.  The Hamiltonian itself can be 
reproduced from the quadratic trace. 

The matrix $F = \sum_{j\not= k} F_{kj}E_{jk}$ , 
which is the residue of $L(z)$ at $z = 0$, stays 
on a coadjoint orbit of $sl(\ell)$ as $t$ varies. 
The phase space of the spin generalization can be 
thereby restricted to the direct product of the 
phase space of Calogero-Moser particles and a 
coadjoint orbit of various dimensions in the 
dual space of $sl(\ell)$.  The lowest dimensional 
non-trivial coadjoint orbit can be parametrized by 
$2\ell$ variables $a_j,b_j$ ($j = 1,\cdots,\ell$) as 
\beqn
    F_{jk} = ig b_j a_k \quad (j \not= k), 
\eeqn
where $g$ is a constant.  These reduced spin 
degrees of freedom, however, can be eliminated 
by a diagonal gauge transformation of the Lax 
equations.  (This does not mean that $a_j$ and 
$b_j$ are non-dynamical.  The elimination 
procedure is done by partially solving the 
equations of motion for those variables.) 
This gauge transformation in turn gives rise 
to non-zero diagonal elements in $M(z)$, and 
the outcome is nothing but the Lax equation 
of the spinless elliptic Calogero-Moser system 
with coupling constant $g$. The spinless system 
is thus embedded in the spin generalization.

\subsubsection{Isomonodromic system} 
There is no substantial difference in the 
construction of an isomonodromic system.  
The equations of motion are given by 
\beqn
    2 \pi i \frac{dq_j}{d\tau} &=& p_j, \quad 
    2 \pi i \frac{dp_j}{d\tau}  =  
        \sum_{k \not= j} \wp'(q_j - q_k) F_{jk} F_{kj}, 
    \nonumber \\
    2 \pi i \frac{dF_{jk}}{d\tau} &=& 
        \sum_{m \not= j} \wp(q_j - q_m) F_{jm} 
      - \sum_{m \not= k} \wp(q_m - q_k) F_{mk}.  
\eeqn
(Terms including $F_{jj}$'s have been eliminated by 
the constraints.) The Lax equation, too, can be written 
in the same form 
\beqn
    2 \pi i \frac{\rd L(z)}{\rd \tau} 
    + \frac{\rd M(z)}{\rd z} 
    = [L(z), M(z)].  
\eeqn
Behind this Lax equation is the ``heat equation'' 
\beqn
    2 \pi i \frac{\rd \sigma(u,z)}{\rd \tau} 
    + \frac{\rd^2 \sigma(u,z)}{\rd u \rd z} 
    = 0 
\eeqn
satisfied by $\sigma(u,z)$. The final piece of 
the ring is the monodromy of $L(z)$ and $M(z)$: 
\beqn
    && L(z + 1) = L(z), \quad M(z + 1) = M(z), 
    \nonumber \\
    && L(z + \tau) = e^{2\pi iQ} L(z) e^{- 2\pi iQ}, 
    \nonumber \\
    && M(z + \tau) = e^{2\pi iQ} \bigl( M(z) 
         + 2 \pi i L(z) \bigr) e^{-2\pi iQ} - 2\pi iP. 
\eeqn
As opposed to the root type Lax pairs, the ordinary 
differential equation 
\beqn
    \frac{dY(z)}{dz} = L(z) Y(z) 
\eeqn
on the torus $E_\tau$ has only one regular singularity 
at $z = 0$. Thus the local monodromy matrix $\Gamma_0$ 
and the global monodromy matrices $\Gamma_\alpha$ and 
$\Gamma_\beta$ are all that are invariant under 
the deformations.

\subsection{Preliminaries for general simple Lie algebra}

Let $\frakg$ be a (complex) simple Lie algebra of 
rank $\ell$, $\frakh$ a Cartan subalgebra, and 
$\Delta$ the associated root system.  The Cartan 
subalgebra induces a root space decomposition of 
$\frakg$: 
\beqn
    \frakg = \frakh \oplus 
      \bigoplus_{\alpha\in\Delta} \frakg_\alpha. 
\eeqn
We choose a basis $\{e_\alpha, h_\mu \mid 
\alpha \in \Delta, \ \mu = 1,\cdots,\ell\}$ 
of $\frakg$ as follows:  
\begin{enumerate}
\item 
$h_\mu$, $\mu = 1,\dots,\ell$, are an orthonormal 
basis of $\frakh$ with respect to the Killing form 
$B: \frakh \times \frakh \to \bbC$, i.e., 
\beqn
    B(h_\mu,h_\nu) = \delta_{\mu\nu}.
\eeqn
The Killing form induces an isomorphism 
$\frakh^* = Hom(\frakh,\bbC) \simeq \frakh$, 
which determines an element $h_\alpha$ for each 
$\alpha \in \frakh^*$. In terms of the basis 
$h_\mu$ of $\frakh$, this map can be written 
explicitly: 
\beqn
    \alpha \mapsto h_\alpha 
    = \sum_{\mu=1}^\ell \alpha(h_\mu) h_\mu, 
\eeqn
\item 
The root subspace $\frakg_\alpha$ is one dimensional. 
$e_\alpha$ is a basis of $\frakg_\alpha$ such that 
\beqn
    [e_\alpha, e_{-\alpha}] = h_\alpha. 
\eeqn
This choice of $e_\alpha$ amounts to the normalization 
\beqn
    B(e_\alpha, e_{-\alpha}) = 1. 
\eeqn
\end{enumerate}
The Lie brackets of the basis elements other than 
$[e_\alpha,e_{-\alpha}]$ now takes the form
\beqn
    [e_\alpha, e_\beta] &=& N_{\alpha,\beta} e_{\alpha+\beta} 
    \quad (\alpha + \beta \not= 0), 
    \nonumber \\ 
    {}[h_\mu, e_\alpha] &=& \alpha(h_\mu) e_\alpha, 
    \nonumber \\ 
    {}[h_\mu, h_\nu] &=& 0. 
\eeqn
The structure constants $N_{\alpha,\beta}$ are 
anti-symmetric with respect to the indices, and 
vanish if $\alpha + \beta \not\in \Delta$.  
The following general relation among the structure 
constants will be used in the course of the proof 
of a Lax equation 

\begin{lemma} 
\beqn
    N_{-\beta,\alpha+\beta} = N_{-\alpha,-\beta} 
    = N_{\alpha+\beta,-\alpha}. 
    \label{eq:three-N}
\eeqn
\end{lemma}

\proof 
If $\alpha = \beta$, this relation is trivially satisfied, 
because all the structure constants vanish.  Let us 
consider the case where $\alpha \not= \beta$.  By the 
Jacobi identity, we have 
\beqn
    [e_{\alpha+\beta}, [e_{-\alpha},e_{-\beta}]] 
    = [[e_{\alpha+\beta},e_{-\alpha}], e_{-\beta}] 
    + [e_{-\alpha}, [e_{\alpha+\beta},e_{-\beta}]]. 
    \nonumber 
\eeqn
This implies the identity 
\beqn
    N_{-\alpha,-\beta} h_{\alpha+\beta} 
    = N_{\alpha+\beta,-\alpha} h_\beta 
    - N_{\alpha+\beta,-\beta} h_\alpha, 
    \nonumber 
\eeqn
which, by the relation 
$h_{\alpha+\beta} = h_\alpha + h_\beta$, 
can be rewritten 
\beqn
      (N_{-\alpha,-\beta} + N_{\alpha+\beta,-\beta}) h_\alpha 
    + (N_{-\alpha,-\beta} - N_{\alpha+\beta,-\alpha}) h_\beta 
    = 0.  
    \nonumber 
\eeqn
Since we have assumed that $\alpha \not= \beta$, 
$h_\alpha$ and $h_\beta$ are linearly independent, 
so that the two coefficients in this linear retion 
are equal to zero. \qed

We can now specify the classical spin variables 
for a general simple Lie algebra. Those spin 
variables, by definition, are coordinates of 
the dual space $\frakg^* = Hom(\frakg,\bbC)$. 
Let $F_\alpha$ and $G_\mu$ be the coordinates 
dual to the above basis $e_\alpha$ and $h_\mu$.  
In other words, they are the coefficients of 
$e_\alpha$ and $h_\mu$ in the linear combination 
\beqn
    \sum_{\alpha\in\Delta} F_{-\alpha} e_\alpha 
    + \sum_{\mu=1}^\ell G_\mu h_\mu 
\eeqn
that realizes the isomorphism $\frakg^* \simeq \frakg$ 
induced by the Killing form.  The Kostant-Kirillov 
Poisson structure on $\frakg^*$ determine the Poisson 
brackets of these spin variables, which take the same 
form as the Lie brackets of the Lie algebra basis: 
\beqn
    \{F_\alpha, F_{-\alpha}\} &=& G_\alpha 
      = \sum_{\mu=1}^\ell \alpha(h_\mu) G_\mu, 
    \nonumber \\
    \{F_\alpha, F_\beta\} &=& N_{\alpha,\beta}F_{\alpha+\beta}  
    \quad (\alpha + \beta \not= 0), 
    \nonumber \\
    \{G_\mu, F_\alpha\} &=& \alpha(h_\mu) F_\alpha, 
    \nonumber \\
    \{G_\mu, G_\nu\} &=& 0. 
\eeqn

\subsection{Spin generalization for general simple Lie algebra} 

\subsubsection{Hamiltonian formalism}
The spin generalization based on $\frakg$, 
too, is a constrained Hamiltonian system 
defined on $\frakh \times \frakh \times \frakg^*$ 
by the Hamiltonian 
\beqn
    \calH = \frac{1}{2} B(p,p) 
      - \frac{1}{2} \sum_{\alpha\in\Delta} 
          \wp(\alpha(q)) F_{-\alpha} F_\alpha 
\eeqn 
and the constraints 
\beqn
    G_\mu = 0 \quad (\mu = 1,\cdots,\ell). 
\eeqn
Here $q$ and $p$ are understood to take values in 
$\frakh$. $B(p,q)$ and $\alpha(q)$ amount to 
$p \cdot p$ and $\alpha \cdot q$ in the models 
based on root systems.  Let us use the same 
``dot notation'' for the Killing form 
$\frakh \times \frakh \to \bbC$ and the pairing 
$\frakh^* \times \frakh \to \bbC$. The Hamiltonian 
then takes a more familiar form: 
\beqn
    \calH = \frac{1}{2} p \cdot p 
      - \frac{1}{2} \sum_{\alpha\in\Delta} 
          \wp(\alpha \cdot q) F_{-\alpha} F_\alpha 
\eeqn

The equations of motion can be readily written 
down in the language of the coordinates 
$q_\mu = q \cdot h_\mu$ and momenta 
$p_\mu = p \cdot h_\mu$ of Calogero-Moser 
particles and the spin variables 
$F_\alpha$ and $G_\mu$ on $\frakg^*$: 
\beqn
    \frac{dq_\mu}{dt} &=& p_\mu, 
    \nonumber \\
    \frac{dp_\mu}{dt} &=& 
      - \frac{1}{2} \sum_{\alpha\in\Delta} 
        \alpha \cdot h_\mu \wp'(\alpha \cdot q) 
        F_{-\alpha} F_\alpha, 
    \nonumber \\
    \frac{dF_\alpha}{dt} &=& 
      - \sum_{\beta\in\Delta,\alpha-\beta\in\Delta} 
        \wp(\beta \cdot q) 
        F_{\alpha-\beta} F_\beta N_{\alpha,-\beta} 
      - \wp(\alpha \cdot q) G_\alpha F_\alpha, 
    \nonumber \\
    \frac{dG_\mu}{dt} &=& 0.      
\eeqn
In particular, the diagonal elements $G_\mu$ of 
the spin variables are conserved quantities. 
One can thereby safely put the aforementioned 
constraints.

\subsubsection{Lax pair} 
The integrability of our spin generalization 
is ensured by the existence of a Lax pair 
as follows. 

\begin{proposition}
Let $V$ be any finite dimensional representation of 
$\frakg$, and $E_\alpha$ and $H_\mu$ the endomorphisms 
on $V$ that represent $e_\alpha$ and $h_\mu$.  
Then the endomorphisms 
\beqn
    L(z) &=& 
    P + \sum_{\alpha\in\Delta} 
        \sigma(\alpha \cdot q, z) F_{-\alpha} E_\alpha, 
    \quad P = \sum_{\mu=1}^\ell p_\mu H_\mu, 
    \nonumber \\
    M(z) &=& 
      - \sum_{\alpha\in\Delta} 
        \sigma(\alpha \cdot q, z) 
           \bigl( \rho(\alpha \cdot q) 
           + \rho(z - \alpha \cdot q) \bigr) 
           F_{-\alpha} E_\alpha 
\eeqn
on $V$ satisfy the Lax equation 
\beqn
    \frac{\rd L(z)}{\rd t} = [L(z), M(z)]. 
\eeqn
\end{proposition}

\proof 
Using the equations of motion and the constraints, 
one can express the $t$-derivative of the $L$-matrix as 
\beqn
    \frac{\rd L(z)}{\rd t} = I + II + III, 
\eeqn
where 
\beqn
    I &=& \sum_{\mu=1}^\ell \frac{dp_\mu}{dt} H_\mu 
      = - \frac{1}{2} \sum_{\alpha\in\Delta} 
            \wp'(\alpha \cdot q) F_{-\alpha} 
            F_\alpha H_\alpha, 
    \nonumber \\
    II &=& \sum_{\alpha\in\Delta} \sum_{\mu=1}^\ell 
           \frac{d\alpha \cdot q}{dt} 
           \left.\frac{\rd \sigma(u,z)}{\rd u}
             \right|_{u=\alpha\cdot q} 
           F_{-\alpha} E_\alpha 
    \nonumber \\
      &=& - \sum_{\alpha\in\Delta} 
            \alpha \cdot \alpha \sigma(\alpha \cdot q, z) 
            \bigl( \rho(\alpha \cdot q) 
              + \rho(z - \alpha \cdot q)\bigr) 
            F_{-\alpha} E_\alpha, 
    \nonumber \\
    III &=& \sum_{\alpha\in\Delta} \sigma(\alpha \cdot q, z) 
            \frac{dF_{-\alpha}}{dt} E_\alpha
    \nonumber \\
      &=& - \sum_{\alpha,\beta\in\Delta,\alpha+\beta\not= 0} 
            \sigma(\alpha \cdot q, z) \wp(\beta \cdot q) 
            F_{-\alpha-\beta} F_\beta 
            N_{-\alpha,-\beta} E_\alpha. 
    \nonumber 
\eeqn
Similarly, the commutator of the Lax pair can be 
written 
\beqn
    [L(z), M(z)] = IV + V + VI, 
\eeqn
where $VI$ stands for terms from the commutator $[P, M(z)]$, 
\beqn
    IV &=& - \sum_{\alpha\in\Delta} \sigma(\alpha \cdot q, z) 
             \bigl( \rho(\alpha \cdot q) 
               + \rho(z - \alpha \cdot q) \bigr) 
             F_{-\alpha} [P, E_\alpha] 
    \nonumber \\
      &=& - \sum_{\alpha\in\Delta} \sigma(\alpha \cdot q, z) 
            \bigl( \rho(\alpha \cdot q) 
               + \rho(z - \alpha \cdot q) \bigr) 
            \alpha \cdot p F_{-\alpha} E_\alpha, 
    \nonumber 
\eeqn
and $V + VI$ are the the other terms grouped into the 
Cartan part ($V$) and the off-Cartan part ($VI$), 
\beqn
    V &=& - \sum_{\alpha\in\Delta} 
            \sigma(-\alpha \cdot q, z) \sigma(\alpha \cdot q, z) 
            \bigl( \rho(\alpha \cdot q) 
              + \rho(z - \alpha \cdot q) \bigr) 
            F_{-\alpha} F_\alpha 
            [E_\alpha, E_{-\alpha}] 
    \nonumber \\
      &=& - \sum_{\alpha\in\Delta} 
            \sigma(-\alpha \cdot q, z) \sigma(\alpha \cdot q, z)  
            \bigl( \rho(\alpha \cdot q) 
              + \rho(z - \alpha \cdot q) \bigr) 
            F_{-\alpha} F_\alpha H_\alpha, 
    \nonumber \\
    VI &=& - \sum_{\alpha,\beta\in\Delta,\alpha+\beta\not= 0} 
             \sigma(\alpha \cdot q, q) \sigma(\beta \cdot q, z)  
             \bigl( \rho(\beta \cdot q) 
               + \rho(z - \beta \cdot q) \bigr) 
             F_{-\alpha} F_\alpha [E_\alpha, E_\beta] 
    \nonumber \\
      &=& - \sum_{\alpha,\beta\in\Delta,\alpha+\beta\not= 0} 
            \sigma(\alpha \cdot q, z) \sigma(\beta \cdot q, z) 
            \bigl( \rho(\beta \cdot q) 
              + \rho(z - \beta \cdot q) \bigr) 
            F_{-\alpha} F_{-\beta} 
            N_{\alpha,\beta} E_{\alpha+\beta}. 
    \nonumber 
\eeqn
It is obvious that $IV = II$. Using 
(\ref{eq:sigma-zero-sum-rule}), 
we can readily see that $V = I$.  
Thus it remains to prove that $VI = III$. 
This is achieved as follows: 
\beqn
    VI &=& - \frac{1}{2} 
             \sum_{\alpha,\beta\in\Delta,\alpha+\beta\not= 0} 
             \sigma(\alpha \cdot q, z) \sigma(\beta \cdot q, z) 
             \Bigl( \rho(\beta \cdot q) 
               + \rho(z - \beta \cdot q) 
    \nonumber \\ 
      &&      - \rho(z - \alpha \cdot q) 
               - \rho(\alpha \cdot q) \Bigr)
             F_{-\alpha} F_{-\beta} 
             N_{\alpha,\beta} E_{\alpha+\beta} 
    \nonumber \\
      && [\mbox{symmetrized with respect to $\alpha$ 
          and $\beta$}]
    \nonumber \\
      &=&  - \frac{1}{2} 
             \sum_{\alpha,\beta\in\Delta,\alpha+\beta\not= 0} 
             \sigma((\alpha+\beta) \cdot q, z) 
             \bigl( \wp(\alpha \cdot q) 
               - \wp(\beta \cdot q) \bigr) 
             F_{-\alpha} F_{-\beta} 
             N_{\alpha,\beta} E_{\alpha+\beta}. 
    \nonumber \\
      && [\mbox{(\ref{eq:sigma-sum-rule}) is used}] 
    \nonumber \\
      &=& \sum_{\alpha,\beta\in\Delta, \alpha+\beta\not= 0} 
             \sigma((\alpha+\beta)\cdot q, z) 
             \wp(\beta \cdot q) 
             F_{-\alpha} F_{-\beta} 
             N_{\alpha,\beta} E_{\alpha+\beta} 
    \nonumber \\
      && [\mbox{asymmetrized with respect to $\alpha$ 
          and $\beta$}] 
    \nonumber \\
      &=& \sum_{\alpha,\beta\in\Delta, \alpha+\beta\not= 0} 
             \sigma(\alpha \cdot q, z) \wp(\beta \cdot q) 
             F_{-\alpha,-\beta} F_\beta 
             N_{\alpha+\beta,-\beta} E_\alpha. 
    \nonumber \\
      && [\mbox{substituting $\beta \to -\beta$ 
          and $\alpha \to \alpha + \beta$}] 
    \nonumber 
\eeqn
Finally using the identity $N_{\alpha+\beta,-\beta} 
= - N_{-\alpha,-\beta}$, cf. (\ref{eq:three-N}), 
we find that the last sum is equal to $III$.  \qed 

Note that the above proof persists to be meaningful 
if $E_\alpha$ and $H_\mu$ are replaced by the Lie 
algebra elements $e_\alpha$ and $h_\mu$.  In other 
words, the Lax equation actually lives in the 
Lie algebra $\frakg$ itself rather than its 
representations.  This resembles the case of the 
Toda systems.

\subsubsection{Isomonodromic System} 
The passage to an isomonodromic analogue is 
straightforward.  Replacing 
$d/dt \to 2 \pi i d/d\tau$, one obtains the 
non-autonomous system 
\beqn
    2 \pi i \frac{dq_\mu}{d\tau} &=& p_\mu, 
    \nonumber \\
    2 \pi i \frac{dp_\mu}{d\tau} &=& 
      - \frac{1}{2} \sum_{\alpha\in\Delta} 
        \alpha \cdot h_\mu \wp'(\alpha \cdot q) 
        F_{-\alpha} F_\alpha, 
    \nonumber \\
    2 \pi i \frac{dF_\alpha}{d\tau} &=& 
      - \sum_{\beta\in\Delta, \alpha-\beta\in\Delta} 
        \wp(\beta \cdot q) 
        F_{\alpha-\beta} F_\beta N_{\alpha,-\beta}.  
\eeqn
(Terms icluding $G_\mu$'s have been eliminated by 
the constraints.)  These equations can be converted 
to the Lax equation
\beqn
    2 \pi i \frac{\rd L(z)}{\rd \tau} 
    + \frac{\rd M(z)}{\rd z}
    = [L(z), M(z)]. 
\eeqn
The monodromy of $L(z)$ and $M(z)$, too, takes 
the same form: 
\beqn
    && L(z + 1) = L(z), \quad M(z + 1) = M(z), 
    \nonumber \\
    && L(z + \tau) = e^{2\pi iQ} L(z) e^{- 2\pi iQ}, 
    \nonumber \\
    && M(z + \tau) = e^{2\pi iQ} \bigl( M(z) 
         + 2\pi i L(z) \bigr) e^{-2\pi iQ} - 2\pi iP, 
\eeqn
where $Q = \sum_{\mu=1}^\ell q_\mu H_\mu$.  
The Lax equation implies that the monodromy data of the 
ordinary differential equation 
\beqn
    \frac{dY(z)}{dz} = L(z) Y(z) 
\eeqn 
on the torus $E_\tau$ is invariant as $\tau$ varies. 
$Y(z)$ now take values in the representation space $V$; 
the monodromy around a singular point or of a cycle of 
$E_\tau$ is represented by a linear transformation on 
$V$.  The ordinary differential equation has a regular 
singularity at $z = 0$ only.  The local monodromy 
around this singular point is a linear transformation 
$\Gamma_0 \in GL(V)$.  Similarly, the global monodromy 
along the $\alpha$ and $\beta$ cycles give 
$\Gamma_\alpha,\Gamma_\beta \in GL(V)$.  These linear 
transformations $\Gamma_0$, $\Gamma_\alpha$ and 
$\Gamma_\beta$ are the monodromy data that are left 
invariant.

%%%%%%%%%%%%%%%%%%%%%%%%%%%%%%%%%%%%%%%%%%%%%%%%%%%%%%%%%%%%
\section{Conclusion} 
\setcounter{equation}{0}

We have thus demonstrated that various models of 
the elliptic Calogero-Moser systems are accompanied 
with an isomonodromic partner.  A technical clue 
is the choice of fundamental functions $x(u,z)$, 
$y(u,z)$, etc. in the Lax pair $L(z)$ and $M(z)$. 
For $L(z)$ and $M(z)$ to give an isomonodromic 
Lax pair, these functions are required to satisfy 
a kind of ``heat equation'' besides the functional 
equations.  We have illustrated the construction 
of the isomonodromic Lax pair for several typical 
cases --- the Lax pair of the $A_{\ell-1}$ mode 
in the vector representation, the root type Lax 
pair for various untwisted and twisted models, and 
the Lax pair of the spin generalizations.  

The most interesting case in the context of Manin's 
equation is the root type Lax pair for the extended 
twisted $BC_\ell$ model (or, equivalently, the 
Inozemtsev system). The root type Lax pair based 
on short roots of the $BC_\ell$ root system consists 
of $2\ell \times 2\ell$ matrices.  

The construction of a Lax pair, however, is merely 
the first step towards a full understanding of Manin's 
equation and its possible generalizations.  The next 
isse is to elucidate the meaning of the affine Weyl 
group symmetries, various special solutions, etc. 
in this framework.  Recent works by Noumi and Yamada 
\cite{bib:noumi-yamada98}, Deift, Its, Kapaev and 
Zhou \cite{bib:deift-etal98} and Kitaev and Korotkin 
\cite{bib:kita-koro98} are very suggestive in this 
respect.  

The spin generalization that we have discussed is 
a special case of a more general multi-spin system, 
i.e., the elliptic Calogero-Moser systems coupled to 
``Gaudin spins'' sitting at the punctures of a 
punctured torus \cite{bib:nekrasov95,bib:enri-rubt95}.  
This is the Hitchin system on a punctured torus; 
we have considered the case with only one puncture 
located at $z = 0$.  It is rather straightforward, 
though more complicated, to generalize our Lax pair 
to the multi-spin generalization.  This gives a 
generalization, to other simple Lie groups, of 
the $SU(2)$ isomonodromic system of Korotkin and 
Samtleben \cite{bib:koro-samt95}. The dynamical 
$r$-matrix in the work of Felder and Wieczerkowski 
\cite{bib:feld-wiec94} plays a central role here. 
We shall report this result elsewhere.

%%%%%%%%%%%%%%%%%%%%%%%%%%%%%%%%%%%%%%%%%%%%%%%%%%%%%%%%%%%%
\section*{Acknowledgements}

I am indebted to Ryu Sasaki for a number of crucial ideas 
on the elliptic Calogero-Moser systems.  I am also grateful 
to Shingo Kawai and Kazuo Okamoto for valuable comments 
on isomonodromic systems on the torus.  Finally, I would 
like to thank Koji Hasegawa, Gen Kuroki, Takashi Takebe 
and Yasuhiko Yamada for discussions on many aspects of 
isospectral and isomonodromic systems.  
This work is partly supported by the Grant-in-Aid for 
Scientific Research (No. 10640165) from the Ministry of 
Education, Science and Culture.

%%%%%%%%%%%%%%%%%%%%%%%%%%%%%%%%%%%%%%%%%%%%%%%%%%%%%%%%%%%%%
\appendix
\renewcommand{\theequation}{\Alph{section}.\arabic{equation}}
%%%%%%%%%%%%%%%%%%%%%%%%%%%%%%%%%%%%%%%%%%%%%%%%%%%%%%%%%%%%%

\section{Proof of Functional Identities and Heat Equation 
for Untwisted Models} 
\setcounter{equation}{0}

\subsection{Proof of (\ref{eq:sum-rule})} 

Let $f(u,v,z)$ denote the difference of both hand sides 
of (\ref{eq:sum-rule}): 
\beqn
    f(u,v,z) = x(u,z) y(v,z) - y(u,z) x(v,z) 
       - x(u+v,z) \bigl( \wp(u) - \wp(v) \bigr). 
\eeqn
This function turns out to have the following analytical 
properties: 
\begin{enumerate}
\item 
$f(u,v,z)$ has the same quasi-periodicity as $x(u,z)$ 
on the $u$ plane, i.e., 
\beqn
    f(u + 1, v, z) = f(u,v,z), \quad 
    f(u + \tau, v, z) = e^{2\pi iz} f(u,v,z). 
\eeqn
\item
$f(u,v,z)$ is an entire function on the $u$ plane. 
\end{enumerate}
The first property is obvious from the quasi-periodicity 
of $x(u,z)$ and the periodicity of $\wp(u)$.  Furthermore, 
poles of $f(u,v,z)$ can appear only at the lattice points 
$u = m + n\tau$ ($m,n \in \bbZ$) on the $u$ plane.  
Therefore, in order to verify the second property, 
we have only to show that $f(u,v,z)$ is non-singular at 
these points. Actually, because of the quasi-periodicity, 
it is sufficient to consider the point $u = 0$ only.  
As $u \to 0$, the singular terms $x(u,z)$, $y(u,z)$ and 
$\wp(u)$ in $f(u,v,z)$ behave as 
\beqn
    x(u,z) &=& \frac{1}{u} + O(1), 
    \nonumber \\
    y(u,z) &=& - \frac{1}{u^2} + O(1), 
    \nonumber \\
    \wp(u) &=& \frac{1}{u^2} + O(u^2) 
\eeqn
so that 
\beqn
    f(u,v,z) 
    &=&  \left(\frac{1}{u} + O(1)\right) y(v,z) 
       - \left(-\frac{1}{u^2} + O(1)\right) x(v,z) 
    \nonumber \\
    && - \left(x(u,z) + y(u,z)u + O(u^2)\right) 
         \left(\frac{1}{u^2} - \wp(v) + O(u^2)\right) 
    \nonumber \\
    &=& O(1). 
\eeqn
We can thus verify the above two properties of 
$f(u,v,z)$. 

Actually, any function with these two properties 
should vanish identically. This can be seen in 
several different ways.  The shortest will be to 
resort to algebraic geometry of line bundles on 
the torus $E_\tau$.  A more elementary proof is 
to consider the quotient $f(u,v,z)/x(u,z)$.  
This quotient is a doubly-periodic meromorphic 
function, and all possible poles are located 
at the lattice points $u = m + n\tau$ 
($m,n \in \bbZ$), and at most of first order.  
In other words, $f(u,v,z)/x(u,z)$ is a meromorphic 
function on the torus with the only possible pole 
at $u = 0$, but the order of pole cannot be greater 
than one. Such a function has to be a constant. 
On the other hand, because of the pole of $x(u,z)$ 
at $u = 0$, $f(u,v,z)/x(u,z)$ has a zero at $u = 0$. 
Therefore the constant should be equal to zero.

\subsection{Proof of (\ref{eq:zero-sum-rule}) and 
(\ref{eq:factor})}

(\ref{eq:zero-sum-rule}) can be readily derived from 
(\ref{eq:sum-rule}) by letting $v \to -u$. 
Let us consider  (\ref{eq:factor}). 
By (\ref{eq:zero-sum-rule}), 
\beqn
    \frac{\rd}{\rd u}\Bigl( x(u,z) x(-u,z) \Bigr) 
    = - x(u,z) y(-u,z) + y(u,z) x(-u,z) 
    = - \wp'(u). 
\eeqn
Consequently, 
\beqn
    x(u,z) x(-u,z) = - \wp(u) 
    + (\mbox{independent of $u$}). 
\eeqn
Since $x(u,z) = - x(z,u) = - x(-u,-z)$, the left hand side 
of the last relation is in fact an anti-symmetric function 
of $u$ and $z$.  Therefore, 
\beqn
    x(u,z) x(-u,z) = \wp(z) - \wp(u) + \const 
\eeqn
Now consider the limit as $u \to z$.  Both $x(u,z)x(-u,z)$ 
and $\wp(z) - \wp(u)$ tend to zero in this limit. Thus 
the constant on the right hand side has to be zero.

\subsection{Proof of (\ref{eq:x-heat-eq}) } 

Let us rewrite the both hand sides of (\ref{eq:x-heat-eq}) 
into a more accessible form. Differentiating $x(u,z)$ by 
$\tau$ gives 
\beqn
    \frac{\rd x(u,z)}{\rd \tau} = 
    x(u,z) \frac{\rd}{\rd \tau} \Bigl( 
      \log\theta_1(z - u) + \log\theta_1'(0) 
      - \log\theta_1(z) - \log\theta_1(u) \Bigr) . 
\eeqn
By the heat equation (\ref{eq:theta-heat-eq}) of the 
Jacobi theta function, 
\beqn
    4 \pi i \frac{\rd}{\rd \tau} \theta_1(u) 
    = \frac{\theta_1''(u)}{\theta_1(u)} 
    = \frac{\rd}{\rd u}\left( 
          \frac{\theta_1'(u)}{\theta_1(u)} \right)
      + \left( \frac{\theta_1'(u)}{\theta_1(u)} \right)^2 
    = \rho'(u) + \rho(u)^2. 
\eeqn
Letting $u \to 0$ and recalling the singular behavior 
of $\rho(u)$ at $u = 0$, we obtain 
\beqn
    4 \pi i \frac{\rd}{\rd \tau} \log\theta_1'(0) 
    = \lim_{u\to 0} \bigl( \rho'(u) + \rho(u) \bigl) 
    = \frac{\theta_1'''(0)}{\theta_1'(0)}. 
\eeqn
Plugging these formulae into the above expression of 
$\rd x(u,z)/\tau$ gives 
\beqn
    4 \pi i \frac{\rd x(u,z)}{\rd \tau} = x(u,z) f(u,z), 
\eeqn
where 
\beqn
    f(u,z) = \rho'(z - u) + \rho(z - u)^2 
      + \frac{\theta_1'''(0)}{\theta_1'(0)}
      - \rho'(z) - \rho(z)^2 - \rho'(u) - \rho(u)^2. 
\eeqn
On the other hand, we have 
\beqn
    \frac{\rd x(u,z)}{\rd u \rd z} 
    = - \frac{\rd}{\rd z} \Bigl( x(u,z) \bigl(\rho(u) 
        + \rho(z - u) \bigr) \Bigr) 
    = - x(u,z) g(u,z), 
\eeqn
where 
\beqn
    g(u,z) = \bigl(\rho(z - u) - \rho(z)\bigr) 
             \bigl(\rho(u) + \rho(z - u)\bigr) 
           + \rho'(z - u). 
\eeqn

The goal is to verify that $f(u,z) = 2 g(u,z)$. 
It is sufficient to prove the following two properties 
of $f(u,z) - 2g(u,z)$, because such a function has 
to be identically zero.  
\begin{enumerate}
\item 
$f(u,z) - 2g(u,z)$ is a doubly-periodic function on 
the $u$ plane with primitive periods $1$ and $\tau$. 
\item 
$f(u,z) - 2g(u,z)$ is an entire function, and has 
a zero at $u = 0$. 
\end{enumerate}
The first property is obvious if one notices the 
following quasi-periodicity of $f(u,z)$ and $g(u,z)$: 
\beqn
    f(u+1,z) = f(u,z), &&
    f(u+\tau,z) 
      = f(u,z) + 4 \pi i \bigl(\rho(u) + \rho(z-u) \bigr), 
    \nonumber \\
    g(u+1,z) = g(u,z), &&
    g(u+\tau,z) 
      = g(u,z) + 2 \pi i \bigl(\rho(u) + \rho(z-u) \bigr). 
\eeqn
Let us check the second property.  Possible poles of 
$f(u,z)$ and $g(u,z)$ are located at the two points 
$u = 0$ and $u = z$ of the fundamental domain of the 
period lattice $\bbZ + \tau\bbZ$.  Again recalling 
the singular behavior of $\rho(u)$ at $u = 0$, one 
can confirm by straightforward calculations that 
\beqn
    f(u,z) = O(u), \quad 
    g(u,z) = O(u)  \quad (u \to 0). 
\eeqn
Thus $f(u,z) - 2g(u,z)$ turns out to be non-singular and 
have a zero at $u = 0$.   Similarly, one can see that 
$f(u,z) - 2g(u,z)$ is non-singular at $u = z$.

%%%%%%%%%%%%%%%%%%%%%%%%%%%%%%%%%%%%%%%%%%%%%%%%%%%%%%%%%%%%%
\section{Verification of Lax Pair for Extended Twisted 
$BC_\ell$ Model} 
\setcounter{equation}{0}

To prove the Lax equation, it is sufficient to derive the 
following three equations: 
\beqn
    \frac{\rd X_a(z)}{\rd t} &=& [P, X_a(z)] 
      \quad (a = 1,2,3), 
    \label{eq:bc-lax1} \\
    \frac{dp \cdot \mu}{dt} &=& 
      [X_1(z) + X_2(z) + X_3(z), 
       Y_1(z) + Y_2(z) + Y_3(z)]_{\mu\mu}, 
    \label{eq:bc-lax2}\\
    0 &=& [X_1(z) + X_2(z) + X_3(z), 
       D + Y_1(z) + Y_2(z) + Y_3(z)]_{\mu\nu} 
       \quad (\mu \not= \nu). 
    \label{eq:bc-lax3}
\eeqn
$\mu$ and $\nu$ run over the set $\Delta_s$ of 
short roots. 

The proof of (\ref{eq:bc-lax1}) is quite easy.  
Let us consider the case of $a = 1$.  
The $t$-derivative of $X_1(z)$ can be written 
\beqn
    \frac{\rd X_1(z)}{\rd t} 
    = ig_m \sum_{\alpha\in\Delta_m} \alpha \cdot p 
        y(\alpha \cdot q, z) E(\alpha).     
\eeqn
Using the commutation relation 
$[P, E(\alpha)] = \alpha \cdot p E(\alpha)$, 
one can readily see that the right hand side 
is equal to $[P,X_1(z)]$.  The other two in 
(\ref{eq:bc-lax1}) can be similarly derived. 

The rest of this appendix is devoted to the 
other two equations (\ref{eq:bc-lax2}) and 
(\ref{eq:bc-lax3}).

\subsection{Proof of (\ref{eq:bc-lax2})} 

We calculate the diagonal elements 
\beqn
   [X_a(z), Y_b(z)]_{\mu\mu}
   = \sum_{\nu\in\Delta_s} 
     \Bigl( X_{a,\mu\nu}(z) Y_{b,\nu\mu}(z) 
     - Y_{b,\mu\nu}(z) X_{a,\nu\mu}(z) \Bigr)  
\eeqn
of the nine commutators one-by-one.  

\subsubsection{Vanishing terms} 
Some part of the matrix elements of $X_a(z)$ and 
$Y_b(z)$ turn out to vanish by the nature of 
the $BC_\ell$ root system: 
\beqn
    && X_{1,\mu,-\mu}(z) = Y_{1,\mu,-\mu}(z) = 0, 
    \label{eq:x1y1-zero} \\
    && X_{2,\mu\nu}(z) = Y_{2,\mu\nu}(z) = 0 
       \quad (\mu \not= - \nu), 
    \label{eq:x2y2-zero} \\
    && X_{3,\mu\nu}(z) = Y_{3,\mu\nu}(z) = 0 
       \quad (\mu \not= - \nu). 
    \label{eq:x3y3-zero} 
\eeqn
The first relation is due to the fact that 
$\mu - (-\mu) = 2\mu $ can never be a middle 
root. The second and third relations are obvious 
if one notices that $\mu - \nu$ is a long root 
(or, equivalently, twice a short root) if and 
only if $\mu = - \nu$. 

In particular, 
\beqn
      [X_1(z), Y_2(z)]_{\mu\mu} 
    = [X_1(z), Y_3(z)]_{\mu\mu} 
    = [X_2(z), Y_1(z)]_{\mu\mu} 
    = [X_3(z), Y_1(z)]_{\mu\mu} 
    = 0. 
\eeqn

\subsubsection{Calculation of $[X_1(z),Y_1(z)]_{\mu\mu}$} 
By definition, 
\beqn
    [X_1(z), Y_1(z)]_{\mu\mu} = 
    - g_m^2 \sum_{\nu\in\Delta_s, \mu-\nu\in\Delta_m} 
      \Bigl( 
        x((\mu-\nu) \cdot q, z) y((\nu-\mu) \cdot q, z) 
    \nonumber \\
      - y((\mu-\nu) \cdot q, z) x((\nu-\mu) \cdot q, z) 
      \Bigr). 
\eeqn
We rewrite this sum to a sum over the middle root 
$\alpha = \mu - \nu$.  Since the middle roots $\alpha$ 
of this form are characterized by the condition that 
$\alpha \cdot \mu = 1$, the right hand side can be 
rewritten 
\beqn
    - g_m^2 \sum_{\alpha\in\Delta_m, \alpha\cdot\mu=1} 
      \Bigl( x(\alpha \cdot q, z) y(-\alpha \cdot q, z) 
      - y(\alpha \cdot q, z) x(-\alpha \cdot q, z) \Bigr). 
      \nonumber 
\eeqn
Actually, the possible values of $\alpha \cdot \mu$ 
are limited to $0$ and $\pm 1$ only.  Therefore 
this sum is equal to 
\beqn
    - \frac{g_m^2}{2} \sum_{\alpha\in\Delta_m} 
      \alpha \cdot \mu \Bigl( 
        x(\alpha \cdot q, z) y(-\alpha \cdot q, z) 
      - y(\alpha \cdot q, z) x(-\alpha \cdot q, z) 
      \Bigr). 
      \nonumber 
\eeqn
(The factor $1/2$ compensates the contributions from 
$\alpha \cdot \mu = 1$ and $\alpha \cdot \mu = -1$.) 
Noting that 
$\alpha \cdot \mu = \{ p \cdot \mu, \alpha \cdot q\}$, 
we can express $[X_1(z), Y_1(z)]$ as a Poisson bracket 
of the form 
\beqn
    [X_1(z), Y_1(z)]_{\mu\mu} 
    = \{ p \cdot \mu, V_{11} \}, 
\eeqn
where 
\beqn
    V_{11} = \frac{g_m^2}{2} \sum_{\alpha\in\Delta_m} 
      x(\alpha \cdot q, z) x(-\alpha \cdot q, z). 
\eeqn

\subsubsection{Contributions of other commutators} 
By (\ref{eq:x2y2-zero}) and (\ref{eq:x3y3-zero}), 
the diagonal elements of the other commutators  
are a sum of just two terms: 
\beqn
    [X_a(z), Y_b(z)]_{\mu\mu} 
    = X_{a,\mu,-\mu} Y_{b,-\mu,\mu} 
    - Y_{b,\mu,-\mu} X_{a,-\mu,\mu}. 
\eeqn

Let us consider the case of $a = 2$ and $b = 2$ 
in some detail. By definition, 
\beqn
    && [X_2(z), Y_2(z)]_{\mu\mu} 
    \nonumber \\
    &=& - \Bigl( g_{l1} x(2\mu \cdot q, z) 
               + g_{l2} x^{(2)}(2\mu \cdot q, z) \Bigr) 
          \Bigl( g_{l1} y(-2\mu \cdot q, z) 
               + g_{l2} y^{(2)}(-2\mu \cdot q, z) \Bigr) 
    \nonumber \\
    &&  + \Bigl( g_{l1} y(2\mu \cdot q, z) 
               + g_{l2} y^{(2)}(2\mu \cdot q, z) \Bigr) 
          \Bigl( g_{l1} x(-2\mu \cdot q, z) 
               + g_{l2} x^{(2)}(-2\mu \cdot q, z) \Bigr) . 
    \nonumber 
\eeqn
Since $\alpha = 2\mu$ is a long root, and long roots 
with non-vanishing inner product with $\mu$ are 
$2\mu$ and $-2\mu$ only, the right hand side can be 
rewritten 
\beqn
    && - \frac{1}{4} \sum_{\alpha\in\Delta_l} 
          \alpha \cdot \mu 
          \Bigl( g_{l1} x(\alpha \cdot q, z) 
               + g_{l2} x^{(2)}(\alpha \cdot q, z) \Bigr) 
          \Bigl( g_{l1} y(-\alpha \cdot q, z) 
               + g_{l2} y^{(2)}(-\alpha \cdot q, z) \Bigr) 
    \nonumber \\
    && +  \frac{1}{4} \sum_{\alpha\in\Delta_l} 
          \alpha \cdot \mu 
          \Bigl( g_{l1} y(\alpha \cdot q, z) 
               + g_{l2} y^{(2)}(\alpha \cdot q, z) \Bigr) 
          \Bigl( g_{l1} x(-\alpha \cdot q, z) 
               + g_{l2} x^{(2)}(-\alpha \cdot q, z) \Bigr). 
    \nonumber 
\eeqn
(The factor $1/4$ compensates the contributions from 
$\alpha \cdot \mu = 2$ and $\alpha \cdot \mu = -2$.) 
We can again cast this into a Poisson bracket: 
\beqn
    [X_2(z), Y_2(z)]_{\mu\mu} = \{p \cdot \mu, V_{22}\}, 
\eeqn
where 
\beqn
    V_{22} = \frac{1}{4} \sum_{\alpha\in\Delta_l} 
        \Bigl( g_{l1} x(\alpha \cdot q, z) 
             + g_{l2} x^{(2)}(\alpha \cdot q, z) \Bigr) 
        \Bigl( g_{l1} x(-\alpha \cdot q, z) 
             + g_{l2} x^{(2)}(-\alpha \cdot q, z) \Bigr). 
    \nonumber \\ 
\eeqn

Similarly, one can obtain 
\beqn
    && [X_2(z), Y_3(z)]_{\mu\mu} = \{p \cdot \mu, V_{23}\}, 
    \quad 
       [X_3(z), Y_2(z)]_{\mu\mu} = \{p \cdot \mu, V_{32}\}, 
    \nonumber \\
    && [X_3(z), Y_3(z)]_{\mu\mu} = \{p \cdot \mu, V_{33}\}, 
\eeqn
where 
\beqn
    V_{23} 
    &=& \frac{1}{2} \sum_{\alpha\in\Delta_s} 
        \Bigl( g_{l1} x(2\alpha \cdot q, z) 
             + g_{l2} x^{(2)}(2\alpha \cdot q, z) \Bigr) 
        \Bigl( g_{s1} x(-\alpha \cdot q, 2z) 
             + g_{s2} x^{(1/2)}(-\alpha \cdot q, 2z) \Bigr), 
    \nonumber \\
    V_{32} 
    &=& \frac{1}{2} \sum_{\alpha\in\Delta_s} 
        \Bigl( g_{s1} x(\alpha \cdot q, 2z) 
             + g_{s2} x^{(1/2)}(\alpha \cdot q, 2z) \Bigr) 
        \Bigl( g_{l1} x(-2\alpha \cdot q, z) 
             + g_{l2} x^{(2)}(-2\alpha \cdot q, z) \Bigr), 
    \nonumber \\
    V_{33} 
    &=& \sum_{\alpha\in\Delta_s} 
        \Bigl( g_{s1} x(\alpha \cdot q, 2z) 
             + g_{s2} x^{(1/2)}(\alpha \cdot q, 2z) \Bigr) 
        \Bigl( g_{s1} x(-\alpha \cdot q, 2z) 
             + g_{s2} x^{(1/2)}(-\alpha \cdot q, 2z) \Bigr). 
    \nonumber \\
\eeqn

Collecting the results of these calculations, 
we find that the right hand side of (\ref{eq:bc-lax2}) 
takes the form of the Poisson bracket $\{p \cdot \mu, V\}$, 
where 
\beqn
    V = V_{11} + V_{22} + V_{23} + V_{32} + V_{33}. 
\eeqn

\subsubsection{Writing $V$ in terms of $\wp$ functions} 
The final step is to rewrite $V$ in terms of the Weierstrass 
$\wp$ functions. For $V_{11}$, this can be done by use of 
(\ref{eq:factor}).  The other parts are due to the following 
functional identities: 
\beqn
    && x^{(1/2)}(u,z) x^{(1/2)}(-u,z) 
      = - \wp^{(1/2)}(u) + \wp^{(1/2)}(\frac{z}{2}), 
    \\
    && x^{(2)}(u,z) x^{(2)}(-u,z) 
      = - \wp^{(2)}(u) + \wp^{(2)}(2z), 
    \\ 
    && x(u, 2z) x^{(1/2)}(-u,2z) + x^{(1/2)}(u,2z) x(-u,2z) 
      = - 2 \wp(u) + \const, 
    \label{eq:twisted-factor-1st} \\
    && x(u,2z) x(-2u,z) + x(2u,z) x(-u,2z) 
      = - \wp(u) + \const, 
    \\
    && x(u,2z) x^{(2)}(-2u,z) + x^{(2)}(2u,z) x(-u,2z) 
      = - \wp(u) + \const, 
    \\
    && x^{(1/2)}(u,2z) x(-2u,z) + x(2u,z) x^{(1/2)}(-u,2z) 
      = - \wp^{(1/2)}(u) + \const, 
    \\
    && x^{(1/2)}(u,2z) x^{(2)}(-2u,z) 
      + x^{(2)}(2u,z) x^{(1/2)}(-u,2z) ]
      = - \wp(u) + \const, 
    \\
    && x(u,z) x^{(2)}(-u,z) + x^{(2)}(u,z) x(-u,z) 
      = - 2\wp^{(2)}(u) + \const 
    \label{eq:twisted-factor-last}  
\eeqn
The first two are substantially the same as 
(\ref{eq:factor}) except that the variables 
and the primitive periods are rescaled. 
``$\const$'' in the other identities stand 
for terms that are independent of $u$,  
thereby negligible in the Poisson bracket 
with $p \cdot \mu$; remember that they 
are not absolute constants, but functions 
of $z$ and $\tau$.  We shall prove these 
identities in Appendix C.  Using these 
functional identities, one can see that 
$V$ is equal to the potential part of the 
Hamiltonian $\calH$,  up to non-dynamical 
terms independent of $p$ and $q$. 

To summarize, we have shown that the sum of 
the $(\mu,\mu)$ elements of the nine commutators 
coincides with the Poisson bracket 
$\{p \cdot \mu, V\}$, which is equal to 
$dp\cdot\mu/dt$ by the equations of motion 
of the model.

\subsection{Proof of (\ref{eq:bc-lax3})} 

The proof can be separated into the cases where 
$\nu = - \mu$ and $\nu \not= \pm\mu$. 

\subsubsection{$\nu = - \mu$} 
The vanishing of the $(\mu,-\mu)$ elements of the 
commutators other than $[X_a(z), D]$ ($a = 1,2,3$) 
and $[X_1(z), D]$ is immediate from 
(\ref{eq:x2y2-zero}) and (\ref{eq:x3y3-zero}). 
$[X_a(z),D]_{\mu,-\mu} $ vanishes because of 
the symmetry $D_{-\mu} = D_\mu$.  As for 
$[X_1(z), Y_1(z)]_{\mu,-\mu}$, we have 
\beqn
    [X_1(z), Y_1(z)]_{\mu,-\mu} 
    &=& - g_m^2 \sum_{\nu\in\Delta_s\setminus\{\pm\mu\}} 
          x((\mu-\nu) \cdot q, z) y((\nu+\mu) \cdot q, z) 
    \nonumber \\
    &&  + g_m^2 \sum_{\nu\in\Delta_s\setminus\{\pm\mu\}} 
          y((\mu-\nu) \cdot q, z) x((\nu+\mu) \cdot q, z). 
\eeqn
Bu substituting $\nu \to -\nu$, the second sum on the 
right hand turns out to be identical to the first sum. 
The two sums thus cancel with each other.

\subsubsection{$\nu \not= \pm\mu$}
The following can be readily seen by using 
(\ref{eq:x2y2-zero}) and (\ref{eq:x3y3-zero}): 
\beqn
  && [X_2(z),D]_{\mu\nu} = [X_3(z),D]_{\mu\nu} = 0, 
     \nonumber \\
  && [X_2(z),Y_2(z)]_{\mu\nu} = [X_2(z),Y_3(z)]_{\mu\nu} 
     = [X_3(z),Y_3(z)]_{\mu\nu} = 0. 
\eeqn
The $(\mu,\nu)$ elements of other commutators can be 
calculated as follows: 
\beqn
    [X_1(z), D]_{\mu\nu} 
    &=& - X_{1,\mu\mu}(z) (D_\mu - D_\nu) 
    \nonumber \\
    &=& g_m x((\mu-\nu) \cdot q, z) 
    \nonumber \\
    && \times \Bigl( 
          g_{s1} \wp(\mu \cdot q) 
        + g_{s2} \wp^{(1/2)}(\mu \cdot q) 
        + g_{l1} \wp(2\mu \cdot q) 
        + g_{l2} \wp^{(2)}(2\mu \cdot q) 
    \nonumber \\
    &&  - g_{s1} \wp(\nu \cdot q) 
        - g_{s2} \wp^{(1/2)}(\nu \cdot q) 
        - g_{l1} \wp(2\nu \cdot q) 
        - g_{l2} \wp^{(2)}(2\nu \cdot q) 
    \nonumber \\
    && + \sum_{\lambda\in\Delta_m,\alpha\cdot\mu=1} 
           \wp(\alpha \cdot q) 
       - \sum_{\alpha\in\Delta_m,\alpha\cdot\nu=1}
           \wp(\alpha \cdot q) \Bigr). 
    \nonumber \\
    \\ {} 
    [X_1(z),Y_1(z)]_{\mu\nu}
    &=& \sum_{\lambda\in\Delta_s}
         \Bigl( X_{1,\mu\lambda}(z) Y_{1,\lambda\nu}(z) 
         - Y_{1,\mu\lambda}(z) X_{1,\lambda\nu}(z) \Bigr) 
    \nonumber \\
    &=& - g_m^2 \sum_{\lambda\in\Delta_s\setminus\{\mu,\nu\}} 
          \Bigl( x((\mu-\lambda) \cdot q, z) 
                 y((\lambda-\nu) \cdot q, z) 
    \nonumber \\
    &&         - y((\mu-\lambda) \cdot q, z) 
                 x((\lambda-\nu) \cdot q, z) \Bigr). 
    \nonumber \\
    \\ {} 
    [X_1(z), Y_2(z)]_{\mu\nu} 
    &=& X_{1,\mu,-\nu}(z) Y_{2,-\nu,\nu}(z) 
       - Y_{2,\mu,-\mu}(z) X_{1,-\mu,\nu}(z) 
    \nonumber \\
    &=& - g_m x( (\mu+\nu) \cdot q, z) 
          \Bigl( g_{l1} y(-2\nu \cdot q, z) 
               + g_{l2} y^{(2)}(-2\nu \cdot q, z) \Bigr) 
    \nonumber \\
    &&  + \Bigl( g_{l1} y(2\mu \cdot q, z) 
               + g_{l2} y^{(2)}(2\mu \cdot q, z) \Bigr) 
          g_m x(-(\mu+\nu) \cdot q, z). 
    \nonumber \\
    \\ {} 
    [X_1(z), Y_3(z)]_{\mu\nu} 
    &=& X_{1,\mu,-\nu}(z) Y_{3,-\nu,\nu}(z) 
       - Y_{3,\mu,-\mu}(z) X_{1,-\mu,\nu}(z) 
    \nonumber \\
    &=& - g_m x((\mu+\nu) \cdot q, z) 
          \Bigl( g_{s1} y(-\nu \cdot q, 2z) 
               + g_{s2} y^{(1/2)}(-\nu \cdot q, 2z) \Bigr) 
    \nonumber \\
    &&  + \Bigl( g_{s1} y(\mu \cdot q, 2z) 
               + g_{s2} y^{(1/2)}(\mu \cdot q, 2z) \Bigr) 
          g_m x(-(\mu+\nu) \cdot q, z). 
    \nonumber \\
    \\ {} 
    [X_2(z), Y_1(z)]_{\mu\nu} 
    &=&  X_{2,\mu,-\mu}(z) Y_{1,-\mu,\nu}(z) 
       - Y_{1,\mu,-\nu}(z) X_{2,-\nu,\nu}(z) 
    \nonumber \\
    &=& - \Bigl( g_{l1} x(2\mu \cdot q, z) 
               + g_{l2} x^{(2)}(2\mu \cdot q, z) \Bigr) 
          g_m y(-(\mu+\nu) \cdot q, z) 
    \nonumber \\
    &&  + g_m y((\mu+\nu) \cdot q, z) 
          \Bigl( g_{l1} x(-2\nu \cdot q, z) 
               + g_{l2} x^{(2)}(-2\nu \cdot q, z) \Bigr). 
    \nonumber \\
    \\ {} 
    [X_3(z), Y_1(z)]_{\mu\nu} 
    &=& X_{3,\mu,-\nu}(z) Y_{1,-\nu,\nu}(z) 
       - Y_{1,\nu,-\nu}(z) X_{3,-\nu,\nu}(z) 
    \nonumber \\
    &=& - 2 \Bigl( g_{s1}x(\mu \cdot q, 2z) 
                 + g_{s2}x^{(1/2)}(\mu \cdot q, 2z) \Bigr) 
          g_m y(-(\mu+\nu) \cdot q, z) 
    \nonumber \\
    &&  + 2 g_m y((\mu+\nu) \cdot q, z) 
          \Bigl( g_{s1} x(-\nu \cdot q, 2z) 
               + g_{s2} x^{(1/2)}(-\nu \cdot q, 2z) \Bigr). 
    \nonumber \\
\eeqn

We now sum up all these quantities, regroup terms 
into those multiplied by the same monomial of coupling 
constants, and show the cancellation in each partial 
sum.  There are six monomials of coupling constants 
that can occur --- $g_m^2$, $g_m g_{l1}$, $g_m g_{l2}$, 
$g_m g_{s1}$ and $g_m g_{s2}$.  

Let us consider the terms multiplied by $g_m^2$.  
This is a sum of the following two quantities: 
\beqn
    I &=& 
       x((\mu-\nu) \cdot q, z) 
       \Bigl( \sum_{\alpha\in\Delta_m,\alpha\cdot\mu=1} 
              \wp(\alpha \cdot q) 
            - \sum_{\alpha\in\Delta_m,\alpha\cdot\nu=1}
              \wp(\alpha \cdot q)  \Bigr) 
    \nonumber \\
    II &=&  
        -  \sum_{\lambda\in\Delta_s\setminus\{\mu,\nu\}} 
           \Bigl( x((\mu-\lambda) \cdot q, z) 
                  y((\lambda-\nu) \cdot q, z) 
    \nonumber \\
    &&          - y((\mu-\lambda) \cdot q, z)  
                  x((\lambda-\nu) \cdot q, z) \Bigr). 
    \nonumber 
\eeqn
By the functional identity (\ref{eq:sum-rule}), 
we can rewrite $II$ into a sum over middle roots: 
\beqn
    II 
    &=& - \sum_{\lambda\in\Delta_s\setminus\{\mu,\nu\}}
      x((\mu-\nu) \cdot q, z) 
      \Bigl( \wp((\mu-\lambda) \cdot q) 
           - \wp((\nu-\lambda) \cdot q) \Bigr) 
    \nonumber \\
    &=& - x( (\mu-\nu) \cdot q, z) 
        \Bigl( \sum_{\alpha\in\Delta_m, \alpha\cdot\mu=1} 
               \wp(\alpha \cdot q) 
             - \sum_{\alpha\in\Delta_m,\alpha\cdot\nu=1}
               \wp(\alpha \cdot q) \Bigr). 
    \nonumber 
\eeqn
Here the sum over $\lambda$ has been converted to a sum 
over $\alpha$ by putting $\alpha = \mu - \lambda$ and 
$\alpha = \nu - \lambda$ in the two $\wp$ function in 
the first line.  Note that $\mu$, $\nu$ and $\lambda$ 
are all orthogonal to each other.  We thus find that 
$I + II = 0$.  

For the other partial sums, we use the following 
functional identities, which we shall prove in 
Appendix C: 
\beqn
    &&   x(2u,z) y(-u-v,z) - y(2u,z) x(-u-v,z) 
       + x(u+v,z) y(-2v,z) \nonumber \\
    && - y(u+v,z) x(-2v,z) 
       - x(u-v,z) \bigl(\wp(2u) - \wp(2v)\bigr) = 0, 
    \label{eq:twisted-sum-rule-1st} \\
    &&   x^{(2)}(2u,z) y(-u-v,z) - y^{(2)}(2u,z) x(-u-v,z) 
       + x(u+v,z) y^{(2)}(-2v,z) \nonumber \\
    && - y(u+v,z) x^{(2)}(-2v,z) 
       - x(u-v,z) \bigl(\wp^{(2)}(2u) - \wp^{(2)}(2v)\bigr) = 0, 
    \\
    &&   2x(u,2z) y(-u-v,z) - y(u,2z) x(-u-v,z) 
       + x(u+v,z) y(-v,2z) \nonumber \\
    && - 2y(u+v,z) x(-v,2z) 
       - x(u-v,z) \bigl(\wp(u) - \wp(v)\bigr) = 0, 
    \\
    && 2x^{(1/2)}(u,2z) y(-u-v,z) - y^{(1/2)}(u,2z) x(-u-v,z) 
       + x(u+v,z) y^{(1/2)}(-v,2z) 
    \nonumber \\
    && - 2y(u+v,z) x^{(1/2)}(-v,2z) 
       - x(u-v,z) \bigl(\wp^{(1/2)}(u) - \wp^{(1/2)}(v)\bigr) 
       = 0. 
    \label{eq:twisted-sum-rule-last}
\eeqn
By these functional identities, we can confirm that all the 
partial sums regrouped by $g_m g_{l1}$, $g_m g_{l2}$, 
$g_m g_{s1}$ and $g_m g_{s2}$, respectively, cancel out.

%%%%%%%%%%%%%%%%%%%%%%%%%%%%%%%%%%%%%%%%%%%%%%%%%%%%%%%%%%%%%
\section{Proof of Functional Identities for Twisted Models}
\setcounter{equation}{0}

We here prove the functional identities that we have 
encountered in Appendix B.  Although the proof is optimized 
to our choice of $x(u,z)$, $x^{(1/2)}(u,z)$ and $x^{(2)}(u,z)$, 
the same method can in principle apply to other solutions 
of the functional equations, such as the functions used 
by D'Hoker and Phong \cite{bib:dhoker-phong98} and 
Bordner and Sasaki \cite{bib:bordner-etal98c}.  

\subsection{Analytical properties of $x^{(1/2)}(u,z)$ and 
$x^{(2)}(u,z)$}

The proof of the identities including $x^{(1/2)}(u,z)$ 
and $x^{(2)}(u,z)$, like the proof in Appendix A, is 
based on the analytical properties of those functions. 

\begin{itemize} 
\item $x^{(1/2)}(u,z)$ has the following analytical 
properties: 

\begin{enumerate}
\item 
$x^{(1/2)}(u,z)$ is a meromorphic function of $u$ and $z$.  
The poles on the $u$ plane and the $z$ plane are located 
at the lattice points $u = m/2 + n\tau$ and $z = m + 2n\tau$ 
($m,n \in \bbZ$). 
\item
$x^{(1/2)}(u,z)$ has the following quasi-periodicity: 
\beqn
    x^{(1/2)}(u + \frac{1}{2}, z) = x^{(1/2)}(u,z), && 
    x^{(1/2)}(u + \tau, z) = e^{2 \pi iz} x^{(1/2)}(u,z), 
    \nonumber \\
    x^{(1/2)}(u, z + 1) = x^{(1/2)}(u, z), && 
    x^{(1/2)}(u, z + 2\tau) = e^{4 \pi iz} x^{(1/2)}(u,z). 
\eeqn
\item 
At the origin of the $u$ and $z$ planes, this function 
exhibits the following singular behavior: 
\beqn
    && x^{(1/2)}(u,z) = \frac{1}{u} - 2 \rho(z \mid 2\tau) 
       + O(u) \quad (u \to 0), 
    \nonumber \\
    && x^{(1/2)}(u,z) = - \frac{2}{z} + 2 \rho(2u \mid 2\tau) 
       + O(z) \quad (z \to 0). 
\eeqn
\end{enumerate}

\item $x^{(2)}(u,z)$ has the following analytical properties: 

\begin{enumerate}
\item 
$x^{(2)}(u,z)$ is a meromorphic function of $u$ and $z$.  
The poles on the $u$ plane and the $z$ plane are located 
at the lattice points $u = 2m + n\tau$ and $z = m + n\tau/2$ 
($m,n \in \bbZ$). 
\item
$x^{(2)}(u,z)$ has the following quasi-periodicity: 
\beqn
    x^{(2)}(u + 2, z) = x^{(2)}(u,z), && 
    x^{(2)}(u + \tau, z) = e^{2 \pi iz} x^{(2)}(u,z), 
    \nonumber \\
    x^{(2)}(u, z + 1) = x^{(2)}(u,z), && 
    x^{(2)}(u, z + \frac{\tau}{2}) = e^{\pi iu} x^{(2)}(u,z). 
\eeqn
\item 
At the origin of the $u$ and $z$ planes, this function 
exhibits the following singular behavior: 
\beqn
    && x^{(2)}(u,z) = \frac{1}{u} 
       - \frac{1}{2}\rho(z \mid \frac{\tau}{2}) 
       + O(u) \quad (u \to 0), 
    \nonumber \\
    && x^{(2)}(u,z) = - \frac{1}{2z} 
       + \frac{1}{2} \rho(\frac{u}{2} \mid \frac{\tau}{2}) 
       + O(z) \quad (z \to 0). 
\eeqn
\end{enumerate}
\end{itemize}

\subsection{Proof of (\ref{eq:twisted-sum-rule-1st}) 
-- (\ref{eq:twisted-sum-rule-last})} 

These four identities can be treated in much 
the same way.  Let us illustrate the proof 
for (\ref{eq:twisted-sum-rule-1st}) only.  
Since the line of the proof is almost the same 
as the proof of (\ref{eq:sum-rule}), we show 
an outline of the proof and leave the details 
to the reader.  

Let $f(u,v,z)$ denote the left hand side of 
(\ref{eq:twisted-sum-rule-1st}):
\beqn
    f(u,v,z) &=& 
         x(2u,z) y(-u-v,z) - y(2u,z) x(-u-v,z) 
       + x(u+v,z) y(-2v,z) 
    \nonumber \\
    && - y(u+v,z) x(-2v,z) 
       - x(u-v,z) \bigl(\wp(2u) - \wp(2v)\bigr). 
\eeqn
Our task is to show the following analytic properties 
of $f(u,v,z)$, which imply that this function is 
identically zero:   
\begin{enumerate}
\item 
$f(u,v,z)$ has the quasi-periodicity as follows: 
\beqn
    f(u + 1, v, z) = f(u,v,z), \quad 
    f(u + \tau, v, z) = e^{2\pi iz} f(u,v,z). 
\eeqn
\item 
$f(u,v,z)$ is an entire function on the $u$ plane.  
\end{enumerate}
The first property is immediate from the 
quasi-periodicity of $x(u,z)$, etc.  Furthermore, 
it is obvious from the definition that 
all possible poles of $f(u,v,z)$ on the 
$u$ plane are limited to the lattice points 
$u = m/2 + n\tau/2$ and $u = -v + m + n\tau$ 
($m,n \in \bbZ$).  In view of the quasi-periodicity, 
therefore, we have only to verify that $f(u,v,z)$ 
is non-singular at $u = 0, 1/2, \tau/2, 1/2 + \tau/2$, 
and $-v$. 

The absence of poles at $u = 0, 1/2$ and $-v$ can be 
verified by straightforward calculations on the basis 
of the singular behavior of $x(u,z)$, $x^{(1/2)}(u,z)$ 
and $x^{(2)}(u,z)$ as $u \to 0$.  

In order to examine the points $u = \tau/2$ and 
$u = 1/2 + \tau/2$, one has to examine the 
singular behavior of $x(2u,z)$ and $y(2u,z)$ as 
$u \to \tau/2, 1/2 + \tau/2$.  This can be worked out 
by combining the quasi-periodicity of $x(u,z)$ and 
$y(u,z)$ and their singular behavior as $u \to 0$: 
\begin{enumerate}
\item 
As $u \to \tau/2$, 
\beqn
    x(2u,z) &=& e^{2\pi iz} x(2u - \tau, z) 
    = e^{2\pi iz} \left(\frac{1}{2u-\tau} 
      + O(1)\right), 
    \nonumber \\
    y(2u,z) &=& 2^{2\pi iz} y(2u - \tau, z) 
    = e^{2\pi iz} \left(- \frac{1}{(2u-\tau)^2} 
      + O(1)\right). 
\eeqn
\item 
As $u \to 1/2 + \tau/2$, 
\beqn
    x(2u,z) &=& e^{2\pi iz} x(2u - 1 - \tau, z) 
    = e^{2\pi iz} \left(\frac{1}{2u-1-\tau} 
      + O(1)\right), 
    \nonumber \\
    y(2u,z) &=& e^{2\pi iz} y(2u -1 - \tau, z) 
    = e^{2\pi iz} \left(- \frac{1}{(2u-1-\tau)^2} 
      + O(1)\right). 
\eeqn
\end{enumerate}
Using these observations, one can confirm the 
absence of poles of $f(u,v,z)$ at $u = \tau/2$ 
and $ 1/2 + \tau/2$ by direct calculations. 

We can thus verify that $f(u,v,z)$ is indeed an 
entire function on the $u$ plane.

\subsection{Proof of (\ref{eq:twisted-factor-1st}) 
-- (\ref{eq:twisted-factor-last})} 

Rather than directly proving these identities, let us 
prove them in a differentiated form.  For illustration, 
we consider the first identity (\ref{eq:twisted-factor-1st}). 
Differentiating this identity by $u$ gives 
\beqn
    && x(u,2z) y^{(1/2)}(-u,2z) - y(u,2z) z^{(1/2)}(-u,2z) 
       + x^{(1/2)}(u,2z) y(-u,2z) 
    \nonumber \\
    && - y^{(1/2)}(u,2z) x(-u,2z) = 2 \wp'(u). 
\eeqn
One can prove it directly, repeating the complex analytic 
reasoning that we have presented in other cases.  An 
alternative way is to take the limit, as $v \to u$, of the 
functional identity 
\beqn
    && x(2u,2z) y^{(1/2)}(-u-v,2z) 
       - y(2u,2z) x^{(1/2)}(-u-v,2z) 
       + x^{(1/2)}(u+v,2z) y(-2v,2z) 
    \nonumber \\
    && - y^{(1/2)}(u+v,2z) x(-2v,2z) 
       - x(u-v,z) \bigl(\wp(2u) - \wp(2v)\bigr) = 0. 
\eeqn
(This yields the above identity upon substituting 
$u \to u/2$ and $v \to v/2$.)  This functional identity 
can be derived by the same method as the proof of 
(\ref{eq:twisted-sum-rule-1st}) -- 
(\ref{eq:twisted-sum-rule-last}).  

Similarly, the third and fifth of 
(\ref{eq:twisted-factor-1st}) -- 
(\ref{eq:twisted-factor-last}) are 
obtained from the following functional identities: 
\beqn
    && 2x(u,2z) y^{(2)}(-u-v,z) 
       - y(u,2z) x^{(2)}(-u-v,z) 
       + x^{(2)}(u+v,z) y(-2v,2z) 
    \nonumber \\
    && - 2y^{(2)}(u+v,z) x(-2v,2z) 
       - x(u-v,z) \bigl(\wp(u) - \wp(v)\bigr) = 0, 
    \\
    && 2x^{(1/2)}(u,2z) y^{(2)}(-u-v,z) 
       - y^{(1/2)}(u,2z) x^{(2)}(-u-v,z) 
       + x^{(2)}(u+v,z) y^{(1/2)}(-2v,2z) 
    \nonumber \\
    && - 2y^{(2)}(u+v,z) x^{(1/2)}(-2v,z) 
       - x(u-v,z) \bigl(\wp(u) - \wp(v)\bigr) = 0. 
\eeqn

The second, forth and sixth of 
(\ref{eq:twisted-factor-1st}) -- 
(\ref{eq:twisted-factor-last}) can be similarly 
derived from the last three of 
(\ref{eq:twisted-sum-rule-1st}) -- 
(\ref{eq:twisted-sum-rule-last}).  
This completes the proof of the functional 
identities.  

We conclude this appendix with a comment on 
the ``$\const$'' terms of these identities.  
In principle, these terms can be determined 
by examining the identities at a special point 
of the $u$ plane.  Let us consider, e.g., 
(\ref{eq:twisted-factor-1st}).  At $u = z$, 
the first term on the left hand side vanishes. 
Evaluating the other terms at this point, 
therefore, one finds that 
\beqn
    \const = 2 \wp(z) - x^{(1/2)}(z,2z) x(-z,2z). 
\eeqn
The same formula can be reproduced by substituting 
$u = -z$.  One can similarly derive an explicit 
expression for the other identities.

%%%%%%%%%%%%%%%%%%%%%%%%%%%%%%%%%%%%%%%%%%%%%%%%%%%%%%%%%%%%%%%%%%%%

\end{document}